# LOCAL SOLVABILITY
# FOR A CLASS OF PARTIAL DIFFERENTIAL OPERATORS
# WITH DOUBLE CHARACTERISTICS

MICHAEL CHRIST AND G. E. KARADZHOV

ABSTRACT. A necessary and sufficient condition for local solvability is presented for the linear partial differential operators $-X^2 - Y^2 + ia(x)[X, Y]$ in $\mathbf{R}^3 = \{(x, y, t)\}$, where $X = \partial_x$, $Y = \partial_y + x^k \partial_t$, and $a \in C^\infty(\mathbf{R}^1)$ is real valued, for each positive integer $k$.

## 1. INTRODUCTION

We say that a linear partial differential operator $L$ is locally solvable at the origin if there exists an open neighborhood $V$ of $0$ such that for every $f \in C_0^\infty(V)$ there exists $u \in \mathcal{D}'(V)$ satisfying $Lu = f$ in $V$.

Elliptic (linear) partial differential operators are always locally solvable, as are many others. The example $\partial_x + i(\partial_y + x\partial_t)$ of Lewy [10] demonstrated that local solvability does not always hold, and subsequently a very satisfactory characterization of local solvability was obtained [1, 14] for differential operators of principal type [1].

Rather little is known, in contrast, for operators having multiple characteristics. For complex *constants* $\alpha$, $L_\alpha = \partial_x^2 + (\partial_y + x\partial_t)^2 + i\alpha\partial_t$ in $\mathbf{R}^3$ is locally solvable if and only if $\alpha \notin \{\pm 1, \pm 3, \pm 5, \dots\}$. These operators are closely related to the original example of Lewy, and arise naturally, along with variants in which $\partial_x$ and $\partial_y + x\partial_t$ are replaced by more general vector fields and $\partial_t$ by their commutator, in connection with the boundary Cauchy-Riemann complex for pseudoconvex domains in $\mathbf{C}^2$ [5]. For the similar family of operators $(\partial_x - ib_1 x^k \partial_t)(\partial_x - ib_2 x^k \partial_t) + iax^{k-1}\partial_t$ in $\mathbf{R}^2$, where $k$ is odd and $a, b_1, b_2$ are real constants, local solvability likewise holds [6] if and only if a certain explicit discrete set of parameters is avoided.

The situation for left-invariant, second-order operators on Heisenberg groups has been analyzed in great detail by Müller and Ricci [12, 13]. These operators depend only on finitely many complex parameters, but the situation is far more complicated.

*Date*: November 2, 1995.
Christ's research was supported in part by NSF grant DMS–9306833 and at MSRI by NSF grant DMS–9022140. Karadzhov was supported in part by a Fulbright Foundation grant.

[1]For the case of pseudodifferential operators of principal type see [9] and the references therein.

For certain subfamilies local solvability is dramatically unstable, depending on Diophantine properties of the coefficients, although the detailed representation theory makes possible a complete analysis. But for certain other subfamilies, including the transversally elliptic operators, local solvability holds for a Zariski open set of coefficients.

Left-invariant operators on Lie groups are a rather natural source of examples for illustrative analysis, but an unnatural end goal from the perspective of the general study of partial differential equations. Extension of the analysis of those unstable families of operators mentioned above to non-left-invariant perturbations is not a realistic goal, but we believe that stable, generically solvable behavior persists for transversally elliptic operators, and indeed for larger classes of operators as well. Moreover nonsolvability should be a very exceptional event, far more exceptional than it is within the context of left-invariant operators on groups, even though it will not be feasible to calculate precisely which coefficients are the exceptional ones.

Consider any two real, smooth vector fields $X, Y$ in $\mathbf{R}^3$ such that $X, Y$ and $[X, Y]$ are linearly independent at 0 and define

$$L = -X^2 - Y^2 + i\alpha[X, Y] \tag{1.1}$$

where $\alpha$ is some $C^\infty$ coefficient. If $\alpha(0) \notin E = \{\pm 1, \pm 3, \pm 5, \dots\}$ then $L$ is subelliptic and hence locally solvable. We conjecture that when $\alpha(0) \in E$, then $L$ is locally solvable at the origin for generic[2] $\alpha(x) - \alpha(0)$.

Our aim in this paper is a preliminary investigation of certain operator families depending on infinitely many parameters, rather than on merely finitely many, in which symmetry is partially broken by the addition of lower order terms. We are at present able only to analyze special situations where separation of variables reduces matters to the analysis of certain eigenvalue problems.[3] Let $X = \partial_x$, $Y = \partial_y + x^k \partial_t$ for some integer $k \geq 1$ and assume $a \in C^\infty$ to be real valued.

**Theorem 1.1.** $L = -X^2 - Y^2 + ia(x)[X, Y]$ *is not locally solvable at the origin if and only if one of the following occurs.*

- $k = 1$, $a(0) \in \{\pm 1, \pm 3, \pm 5, \dots\}$ *and* $a^{(m)}(0) = 0$ *for every* $m \geq 1$.

---

[2]To clarify the appropriate definition of "generic", and thereby to quantify the degree to which the nonsolvable operators are exceptional, is an open problem.

[3]In the present paper, separation of variables reduces matters to eigenvalue problems for certain *ordinary* differential operators. Most of our analysis should be susceptible to generalizations where these are replaced by certain globally elliptic partial differential operators, but we do use repeatedly the fact that all eigenspaces of these ordinary differential operators are one dimensional. This fails to hold for natural generalizations to more than one variable, leading to complications that should not be insurmountable.

- $k > 1$, $a(0) \in \{\pm 1\}$ and $a^{(m)}(0) = 0$ for every $m \geq 1$.

The notation $a^{(m)}$ denotes the derivative of order $m$. The special case where $k = 2$ and $a$ is a constant function has recently been treated by Müller [11], while the case where $k = 1$ and $a$ is constant is well known. We are not aware of any prior work on the case of nonconstant $a$ with $a(0) \in \{\pm 1, \pm 3, \pm 5, \dots\}$.

**Corollary 1.2.** *If $a$ is real analytic and nonconstant then $L$ is locally solvable at the origin.*

In the theory for operators of principal type, what is proved is not merely local solvability, but the stronger property of local solvability in $L^2$, which means that for any $f \in L^2$ there exists a locally square integrable function $u$ satisfying $Lu = f$ in a neighborhood of 0. Many of the operators $L$ whose local solvable is asserted here are not locally solvable in $L^2$; see the final remark in section 6. Thus we are often in the more delicate situation where derivatives are lost.

It is tempting to interpret our results as supporting the thesis that within the class of operators $X^2 + Y^2 + i\alpha[X, Y]$, nonsolvability is not only a rather rare phenomenon, but occurs only in situations that are either highly symmetric, or more generally are nearly reducible to a short list of highly symmetric examples by such operations as conjugation with elliptic Fourier integral operators. The fact that the conditions for nonsolvability are less restrictive for the more symmetric case $k = 1$ than for $k > 1$ is consistent with this thesis. However another class of examples analyzed in [3] demonstrates that the situation is subtler than our results might suggest.

A related family of examples may be defined by taking $a$ to be constant but replacing $Y$ by $\partial_y + b(x)\partial_t$. The same method should apply and should yield similar results, but this has not been investigated in detail.

## 2. Preliminaries

Define the partial Fourier transform with respect to $(y, t)$ by
$$\hat{f}(x, \eta, \tau) = \int \int f(x, y, t) e^{i(y\eta + t\tau)} \, d\eta \, d\tau.$$
For each $s \in \mathbf{R}$ define $\Lambda^s$ by
$$(\widehat{\Lambda^s f})(x, \eta, \tau) = (1 + \tau^2)^{s/2} \hat{f}(x, \eta, \tau).$$

Denote by $L^*$ the transpose of $L$; this is the operator obtained by replacing $a$ by $-a$ in the definition of $L$ and hence is unitarily equivalent to $L$ under the transformation $(x, y, t) \mapsto (x, y, -t)$. Local solvability of $L$ at 0 would follow from an inequality
$$\|\Lambda^s L^* \psi\| \geq c \|\psi\| \text{ for all } \psi \in C_0^\infty(U) \tag{2.1}$$

for some $s < \infty$, some neighborhood $U$ of the origin, and some constant $c > 0$, where $\|\cdot\|$ denotes the $L^2$ norm. Indeed, this implies that

$$\|\psi\|_{L^2} \leq c^{-1}\|L^*\psi\|_{H^s} \quad \text{for all } \psi \in C_0^\infty(U),$$

from which it follows by a straightforward duality argument and the Hahn-Banach theorem that for each $f \in L^2$ supported in $U$, there exists $u \in H^{-s}(U)$ satisfying $Lu = f$.

Define the ordinary differential operators

$$A = A_{\eta,\tau} = -\partial_x^2 + (\eta + x^k\tau)^2 - kx^{k-1}a(x)\tau, \tag{2.2}$$

Then $(\widehat{Lf}) \equiv A_{\eta,\tau}\hat{f}$.

The bulk of our analysis is devoted to the proof of an estimate

$$(1+\tau^2)^s\|A_{\eta,\tau}\phi\| \geq c\|\phi\|, \text{ for all } \phi \in C_0^\infty(U) \tag{2.3}$$

for some $s < \infty$, $c > 0$, and small neighborhood $U$ of 0, for all $(\eta, \tau)$ outside of an exceptional set $\Sigma$ having finite Lebesgue measure. Then if $U$ is chosen to have sufficiently small diameter, (2.1) will follow from the following version of the uncertainty principle.

**Lemma 2.1.** *For each $n \geq 1$ there exists $C < \infty$ such that for each $\delta > 0$, for every measurable set $E \subset \mathbf{R}^n$ having Lebesgue measure less than $C^{-1}\delta^{-1}$ and each function $f \in L^2(\mathbf{R}^n)$ supported on a set of measure less than $\delta$,*

$$\|f\|_{L^2} \leq C\|\hat{f}\|_{L^2(\mathbf{R}^n\setminus E)}.$$

*Proof.* We have

$$\|f\|_{L^2}^2 = c\|\hat{f}\|_{L^2(E)}^2 + c\|\hat{f}\|_{L^2(\mathbf{R}^n\setminus E)}^2$$
$$\leq c|E|\cdot\|\hat{f}\|_{L^\infty}^2 + c\|\hat{f}\|_{L^2(\mathbf{R}^n\setminus E)}^2$$
$$\leq c|E|\cdot\|f\|_{L^1}^2 + c\|\hat{f}\|_{L^2(\mathbf{R}^n\setminus E)}^2$$
$$\leq c|E|\delta\|f\|_{L^2}^2 + c\|\hat{f}\|_{L^2(\mathbf{R}^n\setminus E)}^2$$
$$\leq \frac{1}{2}\|f\|_{L^2}^2 + c\|\hat{f}\|_{L^2(\mathbf{R}^n\setminus E)}^2,$$

provided that $C$ is chosen to be sufficiently large. □

Let $\tau_0, \gamma_0, \gamma_1, \delta_0$ be positive constants to be chosen later in the course of the proof; $\tau_0, \gamma_0$ will be large while $\gamma_1, \delta_0$ will be small. For each $k \geq 1$ we decompose $\mathbf{R}^2 =$

$\mathcal{B}_k \cup \mathcal{C}_k \cup \mathcal{D}_k$ where

$$\mathcal{B}_k = \{(\eta,\tau) \in \mathbf{R}^2 : |\tau| \leq \tau_0 \text{ or } |\eta| \geq \gamma_1|\tau|\}$$
$$\mathcal{C}_k = \{(\eta,\tau) \in \mathbf{R}^2 : |\tau| \geq \tau_0 \text{ and } \gamma_0|\tau|^{1/(k+1)} \leq |\eta| \leq \gamma_1|\tau|\}$$
$$\mathcal{D}_k = \{(\eta,\tau) \in \mathbf{R}^2 : |\tau| \geq \tau_0 \text{ and } |\eta| \leq \gamma_0|\tau|^{1/(k+1)}\}.$$

Fix a cutoff function $\zeta \in C_0^\infty(\mathbf{R})$ that is identically equal to one for $|x| \leq \delta_0$ but is supported where $|x| \leq 2\delta_0$. Then when acting on functions supported in $\{|x| \leq \delta_0\}$, $A_{\eta,\tau}$ may equivalently be written as

$$A_{\eta,\tau} = -\partial_s^2 + (\eta + \tau s^k)^2 - [a(0) + b(s)]\tau k s^{k-1} \tag{2.4}$$

where

$$b(s) = (a(s) - a(0))\zeta(s).$$

Henceforth $A_{\eta,\tau}$ denotes always this modified operator.

Although it suffices to prove (2.3) only for functions $\phi$ supported in $\{|x| \leq \delta_0\}$ for a fixed small constant $\delta_0$, it will nonetheless often be useful to regard $A_{\eta,\tau}$ as an operator defined on $L^2(\mathbf{R})$. It is an unbounded operator of Schrödinger type $-\partial_x^2 + V(x)$, where the potential $V$ is continuous and real valued, and $V(x) \to +\infty$ as $|x| \to \infty$ for all $\tau \neq 0$. Thus (disregarding the case $\tau = 0$ as we may since this is a set of parameters of measure zero in $\mathbf{R}^2$) $A_{\eta,\tau}$ is essentially selfadjoint, and has a discrete sequence $\{\mu_j(\eta,\tau)\}$ of eigenvalues, with $\mu_0 < \mu_1 < \cdots \to +\infty$. For $f \in C_0^2(\mathbf{R})$,

$$\|A_{\eta,\tau}f\| \geq \min_j |\mu_j(\eta,\tau)|\, \|f\|, \tag{2.5}$$

so that obtaining a lower bound for $A_{\eta,\tau}$ is equivalent to deriving a lower bound for $\min_j |\mu_j(\eta,\tau)|$. Throughout the paper the symbol $\|\cdot\|$ with no subscript will denote the norm in $L^2(\mathbf{R})$.

For each $\mu \in \mathbf{R}$, the linear space of all solutions of $A_{\eta,\tau}\phi = \mu\phi$ is two-dimensional, but the behavior of $V$ implies [4] the existence of a solution satisfying $|\phi(x)| \to \infty$ as $x \to +\infty$. Therefore each eigenspace of $A_{\eta,\tau}$ is one dimensional.

**Lemma 2.2.** *For each $k$ and each $a \in C^\infty$, given any constants $\tau_0, \gamma_0, \gamma_1 \in \mathbf{R}^+$, there exist $\delta_0 > 0$ and $C < \infty$ such that for all $f \in C_0^2(\mathbf{R})$ supported in $\{|x| \leq \delta_0\}$ and every $(\eta,\tau) \in \mathcal{B}_k$,*

$$\|f\|^2 \leq C\langle A_{\eta,\tau}f, f\rangle. \tag{2.6}$$

*If $k \geq 2$ is even then for any $a \in C^\infty$ and any $\tau_0, \gamma_1 \in \mathbf{R}^+$, $\gamma_0 \in \mathbf{R}$ may be chosen so that for any finite $\delta_0$, the same inequality holds uniformly for all $(\eta,\tau) \in \mathcal{C}_k$ satisfying $\eta \cdot \tau > 0$.*

*Proof.* One has

$$\langle A_{\eta,\tau} f, f \rangle \geq \|f'\|^2 + \int \left[ (\eta + \tau x^k)^2 - C|\tau||x|^{k-1} \right] |f|^2 \, dx$$
$$\geq \int v(x)|f|^2 \, dx$$

for all $f \in C_0^2$ supported where $|x| \leq \delta_0$, where

$$v(x) = 2^{-1}\delta_0^{-2} + (\eta + \tau x^k)^2 - C|\tau||x|^{k-1},$$

because $\|f'\|^2 \geq 2^{-1}\delta_0^{-2}\|f\|^2$ for all such $f$. If $|\tau| \leq \tau_0$ then for all $|x| \leq \delta_0$

$$v(x) \geq 2^{-1}\delta_0^{-2} - C\tau_0\delta_0^{k-1},$$

which is strictly positive provided that $\delta_0$ is chosen to be sufficiently small relative to $\tau_0^{-1}$.

Consider next the case where $|\eta| \geq \gamma_1|\tau|$ and $|\tau| \geq \tau_0$. Then

$$v(x) \geq (\gamma_1 - \delta_0^k)^2 \tau^2 - C\delta_0^{k-1}|\tau|.$$

Given any $\gamma_1, \tau_0 \in \mathbf{R}^+$, $\delta_0$ may be chosen so that this quantity is bounded below by a small constant times $\tau^2$, for all $|\tau| \geq \tau_0$.

Lastly, if $\eta \cdot \tau > 0$ and $k$ is even then

$$v(x) \geq \eta^2 + 2\eta\tau x^k + \tau^2 x^{2k} - C|\tau||x|^{k-1}$$
$$\geq \eta^2 + x^{2k}\tau^2 - C|\tau||x|^{k-1}.$$

Given $C < \infty$ there exists $C' < \infty$ such that

$$C|\tau x^{k-1}| \leq \left(|\tau|^{1/k}|x|\right)^{2k} + C'\left(|\tau|^{1/k}\right)^{2k/(k+1)}$$
$$\leq \tau^2 x^{2k} + C'|\tau|^{2/(k+1)}.$$

If $(\eta, \tau) \in \mathcal{C}_k$ then $|\eta| \geq \gamma_0|\tau|^{1/(k+1)}$. Therefore

$$v(x) \geq \gamma_0^2|\tau|^{2/(k+1)} + x^{2k}\tau^2 - \tau^2 x^{2k} - C'|\tau|^{2/(k+1)}$$
$$\geq \frac{1}{2}\gamma_0^2|\tau|^{2/(k+1)}$$

provided that $\gamma_0$ is chosen to be sufficiently large relative to $\|b\|_{L^\infty}$. $\square$

## 3. Analysis for $(\eta, \tau) \in \mathcal{C}_k$

Recall that for each $k \geq 1$, $(\eta, \tau) \in \mathcal{C}_k$ if $|\tau| \geq \tau_0$ and $\gamma_0 |\tau|^{1/(k+1)} \leq |\eta| \leq \gamma_1 |\tau|$. This section is devoted to the proof of the following result.

**Proposition 3.1.** *Suppose that $a(0) \notin \{\pm 1, \pm 3, \pm 5 \ldots\}$. Then there exist large constants $\gamma_0, \tau_0$ and small constants $\gamma_1, \delta_0, \delta$ such that for all $(\eta, \tau) \in \mathcal{C}_k$,*

$$\|A_{\eta,\tau} f\| \geq \delta |\tau|^{2/(k+1)} \|f\| \tag{3.1}$$

*for all $f \in C_0^2(\mathbf{R})$.*

*Suppose that $a(0) \in \{\pm 1, \pm 3, \pm 5 \ldots\}$. Then there exist large constants $\gamma_0, \tau_0, M$, small constants $\gamma_1, \delta_0$, and a set $E \subset \mathcal{C}_k$ of finite measure such that*

$$\|A_{\eta,\tau} f\| \geq |\tau|^{-M} \|f\| \text{ for all } (\eta, \tau) \in \mathcal{C}_k \backslash E \tag{3.2}$$

*for every $f \in C_0^2(\mathbf{R})$.*

We will prove this for even $k > 1$, for $\{(\eta, \tau) \in \mathcal{C}_k : \tau > 0 \text{ and } \eta < 0\}$, then comment on the changes needed for the other cases. Change variables $(\eta, \tau) \mapsto (z, \epsilon)$ where

$$z^k = -\eta/\tau, \qquad \epsilon^2 = k^{-1} |\eta|^{-(k+1)/k} |\tau|^{1/k},$$

taking the unique positive solution $z$ of the first equation. Then

$$\gamma_0^{1/k} |\tau|^{-1/(k+1)} \leq |z| \leq \gamma_1^{1/k}$$

and

$$k^{-1} \gamma_1^{-(k+1)/k} |\tau|^{-1} \leq \epsilon^2 \leq k^{-1} \gamma_0^{-(k+1)/k};$$

in particular, both $z, \epsilon$ are arbitrarily close to $0$ provided that $\gamma_1, \gamma_0$ are chosen to be sufficiently small and large, respectively. The inverse relations are

$$\eta = -k^{-1} \epsilon^{-2} z^{-1}, \qquad \tau = k^{-1} \epsilon^{-2} z^{-k-1}.$$

Writing $A_{\eta,\tau} = -\partial_s^2 + (\eta + s^k \tau)^2 - k\tau s^{k-1} a(0) - k\tau s^{k-1} b(s)$ and substituting $s = \epsilon z x$ yields $(\epsilon z)^2 A_{\eta,\tau} = B_{z,\epsilon}$ where

$$B_{z,\epsilon} = -\partial_x^2 + q_\epsilon^2(x) - [a(0) + b(\epsilon z x)] q_\epsilon'(x) \tag{3.3}$$

with

$$q_\epsilon(x) = (\epsilon k)^{-1} \left( (\epsilon x)^k - 1 \right).$$

$q_\epsilon(x) = 0$ if and only if $x = \pm \epsilon^{-1}$; $-\partial_s^2 + q_\epsilon^2$ is for small $\epsilon$ a Schrödinger operator whose potential has a double well. In the next subsection it will be shown that for small $\epsilon$, the eigenfunctions of $B_{z,\epsilon}$ corresponding to small (in absolute value) eigenvalues are localized near the wells, and in fact near one well only. The precise behavior of the small eigenvalues of $B_{z,\epsilon}$ is then determined by the Taylor expansions of $q_\epsilon, q'_\epsilon, b(\epsilon z x)$ about a zero of $q_\epsilon$. This effect will be analyzed in a later subsection.

3.1. **Localization estimates and small eigenvalues.** Define

$$w_\epsilon(x) = \begin{cases} |x - \epsilon^{-1}| & \text{if } x \geq 0 \\ |x + \epsilon^{-1}| & \text{if } x \leq 0. \end{cases}$$

Since $B_{z,\epsilon}$ has real coefficients, all its eigenfunctions may be taken to be real valued.

**Lemma 3.2.** *For any $C_1 > 0$ there exist $r > 0$, $C < \infty$ such that for all sufficiently small $(z, \epsilon)$, for any eigenvalue $\lambda \in [-C_1, C_1]$ of $B_{z,\epsilon}$ and any associated eigenfunction $\phi \in L^2(\mathbf{R})$,*

$$\int_{\mathbf{R}} [\phi^2(x) + (\phi')^2(x)] \, e^{r w_\epsilon(x)} \, dx \leq C \|\phi\|^2.$$

*Proof.* Let $N$ be a large positive number and set $w_{\epsilon,N}(x) = \min(N, w_\epsilon(x))$. It suffices to prove the desired estimate with $w_\epsilon$ replaced by $w_{\epsilon,N}$, with $C$ independent of $N$.

Fix $h \in C_0^2(\mathbf{R})$, real valued and identically equal to 1 in some neighborhood of 0. Let $M \in \mathbf{R}^+$ be a second large constant which, like $N$, will eventually tend to $\infty$, assume that $B_{z,\epsilon} \phi = \lambda \phi$ with $|\lambda| \leq C_1$ and consider

$$\int \left( q_\epsilon^2 - a(0) q'_\epsilon - b(\epsilon z x) q'_\epsilon - \lambda \right) \phi^2(x) h^2(M^{-1} x) e^{r w_{\epsilon,N}(x)} \, dx$$

$$= \int \partial_x^2 \phi \cdot \phi e^{r w_{\epsilon,N}} h^2(M^{-1} x) \, dx$$

$$= -\int (\phi')^2 h^2(M^{-1} x) e^{r w_{\epsilon,N}} \, dx$$

$$- r \int \phi \phi' h^2(M^{-1} x) w'_{\epsilon,N} e^{r w_{\epsilon,N}} \, dx$$

$$- 2 M^{-1} \int \phi h'(M^{-1} x) \cdot \phi' h(M^{-1} x) e^{r w_{\epsilon,N}} \, dx.$$

There exists $C_2 < \infty$ such that for all sufficiently small $\epsilon$,

$$q_\epsilon^2(x) \geq (|a(0)| + \|b\|_{L^\infty}) |q'_\epsilon(x)| + C_1 + 1$$

for all $x$ satisfying $\min(|x - \epsilon^{-1}|, |x + \epsilon^{-1}|) \geq C_2$. Since $w_{\epsilon,N}$ is a Lipschitz function satisfying $|w'_{\epsilon,N}| \leq 1$ almost everywhere, applying the Cauchy-Schwarz inequality to the last displayed inequality and moving one term from the right-hand side to the

left yields

$$\int (\phi^2 + (\phi')^2) h^2 (M^{-1}x) e^{rw_{\epsilon,N}} \, dx$$
$$\leq r \int (\phi^2 + (\phi')^2) h^2 (M^{-1}x) e^{rw_{\epsilon,N}} \, dx$$
$$+ M^{-1} \int (\phi')^2 h^2 (M^{-1}x) e^{rw_{\epsilon,N}} \, dx$$
$$+ M^{-1} \int \phi^2 (h')^2 (M^{-1}x) e^{rw_{\epsilon,N}} \, dx$$
$$+ C \int_{|x \pm \epsilon^{-1}| \leq C_2} \phi^2 e^{rw_{\epsilon,N}} \, dx.$$

The first two terms on the right-hand side may be absorbed into the left if $r, M^{-1}$ are sufficiently small. Fixing $\epsilon, N$, the third term on the right is $O(M^{-1})$, so letting $M \to \infty$ yields

$$\int (\phi^2 + (\phi')^2) e^{rw_{\epsilon,N}} \, dx \leq C \int_{|x \pm \epsilon^{-1}| \leq C_2} \phi^2 \, dx,$$

using the fact that $w_{\epsilon,N}(x) = O(1)$ for $|x \pm \epsilon^{-1}| \leq C_2$, uniformly in $\epsilon, N$. Letting $N \to \infty$ and invoking the Lebesgue monotone convergence theorem now yields the conclusion desired. □

When both $\eta, \tau$ are nonzero, $B_{z,\epsilon}$ is a positive scalar multiple of $A_{\eta,\tau}$ and hence is essentially selfadjoint with a discrete sequence of eigenvalues tending to $+\infty$, associated to one dimensional eigenspaces.

**Lemma 3.3.** *If $a(0) \notin \{\pm 1, \pm 3, \pm 5 \ldots\}$ then there exists $\theta > 0$ such that for all $(z, \epsilon) \in \mathbf{R}^2$ satisfying $|(z, \epsilon)| \leq \theta$ and $\epsilon \neq 0$ and for all $f \in C_0^2$,*

$$\|B_{z,\epsilon} f\| \geq \theta \|f\|.$$

*Proof.* Fix $h \in C_0^2(\mathbf{R})$ supported in $[-1, 1]$ and identically equal to 1 on $[-1/2, 1/2]$. Let $\theta < \min(|a(0) \pm 1|, |a(0) \pm 3|, \ldots)$ be a small constant to be chosen below, and consider any small $(z, \epsilon)$. If there exists $f \in C_0^2$ satisfying $\|B_{z,\epsilon} f\| < \theta \|f\|$ then there exist $\lambda \in [-\theta, \theta]$ and $\phi \in L^2$ such that $B_{z,\epsilon} \phi = \lambda \phi$ and $\|\phi\| = 1$. Set

$$\psi(x) = \phi(x) h(\epsilon^{1/2}(x - \epsilon^{-1})) + \phi(x) h(\epsilon^{1/2}(x + \epsilon^{-1})) = \psi^+(x) + \psi^-(x).$$

By Lemma 3.2,

$$\|(B_{z,\epsilon} - \lambda)\psi\| + \|\phi - \psi\| = O(\exp(-c\epsilon^{-1/2}))$$

for some $c > 0$, uniformly in $z$.

For $|x - \epsilon^{-1}| \leq \epsilon^{-1/2}$,

$$q_\epsilon^2(x) - [a(0) + b(\epsilon z x)]q_\epsilon'(x) = (x - \epsilon^{-1})^2 - [a(0) + b(z)] + O(\epsilon^{1/2}|x - \epsilon^{-1}|^2 + \epsilon^{1/2}).$$

This, like all estimates below, holds uniformly for all $|(z, \epsilon)| \leq 1$. Thus

$$\begin{aligned}
\|(B_{z,\epsilon} - \lambda)\psi^+\| \\
= &\|\left(-\partial_x^2 + (x - \epsilon^{-1})^2 - [a(0) + b(z)] - \lambda\right)\psi^+\| \\
&+ O(\epsilon^{1/2})\|((x - \epsilon^{-1})^2 + 1)\psi^+\| \\
= &\|\left(-\partial_x^2 + (x - \epsilon^{-1})^2 - [a(0) + b(z)] - \lambda\right)\psi^+\| + O(\epsilon^{1/2}),
\end{aligned}$$

by Lemma 3.2. Now $-\partial_x^2 + (x - \epsilon^{-1})^2 - [a(0) + b(z)] - \lambda$ has spectrum $\{1, 3, 5, \dots\} - [a(0) + b(z) + \lambda]$. Since $b(0) = 0$, $\theta$ may be chosen to be so small that for all sufficiently small $|z|$, the intersection of this spectrum with $[-2\theta, 2\theta]$ is empty. Equivalently, $\|(-\partial_x^2 + (x - \epsilon^{-1})^2 - [a(0) + b(z)] - \lambda)g\| \geq 2\theta\|g\|$ for all $g \in C_0^2$. Consequently $\|(B_{z,\epsilon} - \lambda)\psi^+\| \geq 2\theta\|\psi^+\| - C\epsilon^{1/2}$.

The same analysis applies to $\psi^-$, with one algebraic change: $q_\epsilon'(-\epsilon^{-1}) = -1$, so the quantity $-a(0) - b(z) - \lambda$ is replaced by $+a(0) + b(z) - \lambda$. Thus $\|(B_{z,\epsilon} - \lambda)\psi^-\| \geq 2\theta\|\psi^-\| - C\epsilon^{1/2}$. Since $\psi^+, \psi^-$ have disjoint supports and $\|\psi^+\|^2 + \|\psi^-\|^2 = 1 + O(\epsilon^{1/2})$, altogether

$$\|(B_{z,\epsilon} - \lambda)\phi\| \geq \|(B_{z,\epsilon} - \lambda)\psi\| - C\epsilon \geq 2\theta\|\psi\| - C\epsilon^{1/2}.$$

If $\epsilon$ is sufficiently small this last quantity is strictly greater than $\theta\|\psi\|$, a contradiction. □

**Lemma 3.4.** *If $a(0) \in \{\pm 1, \pm 3, \pm 5, \dots\}$ then there exists $\theta > 0$ such that for all sufficiently small $(z, \epsilon)$ with $\epsilon \neq 0$, $B_{z,\epsilon}$ has exactly one eigenvalue in $[-\theta, \theta]$ and no eigenvalues satisfying $\theta < |\lambda| \leq 4\theta$.*

*Proof.* Since the change of variables $(x, y, t) \mapsto (x, y, -t)$ has the effect of replacing $a(x)$ by $-a(x)$, it is no loss of generality to assume that $a(0) \geq 0$, so that in the present situation $a(0) \in \{1, 3, 5, \dots\}$. Likewise the case $\epsilon < 0$ reduces to $\epsilon > 0$ by replacing $x$ by $-x$.

Let $\theta > 0$ be a small constant to be chosen below, fix $(z, \epsilon)$, and assume $B_{z,\epsilon}$ to have least two eigenvalues $\lambda_1, \lambda_2 \in [-4\theta, 4\theta]$. Let $\phi_1, \phi_2$ be associated eigenfunctions of norm 1. As in the proof of Lemma 3.3 decompose $\phi_j = \psi_j^+ + \psi_j^- + (\phi_j - \psi_j^+ - \psi_j^-)$. Lemma 3.2 guarantees that $\|\phi_j - \psi_j^+ - \psi_j^-\| = O(\exp(-c\epsilon^{-1/2}))$ for some $c > 0$. Since $q_\epsilon'(-\epsilon^{-1}) = -1$ and $a(0) \geq 1$, the distance from $-a(0) - b(z) + \lambda$ to the spectrum $\{1, 3, 5, \dots\}$ of $-\partial_x^2 + (x + \epsilon^{-1})^2$ is at least $2 - |b(z)| - |\lambda| \geq 1$ for all sufficiently small $z$, assuming that $|\theta| \leq 1/8$. Thus as in the proof of Lemma 3.3, $\|(B_{z,\epsilon} - \lambda_j)\psi_j^-\| \geq$ $c\|\psi_j^-\|$. Now on the support of $\psi_j^-$, $(B_{z,\epsilon} - \lambda_j)\psi_j^- = (B_{z,\epsilon} - \lambda_j)\phi_j + O(\exp(-c\epsilon^{-1/2}))$, by the decay estimate of Lemma 3.2. Since $(B_{z,\epsilon} - \lambda_j)\phi_j \equiv 0$, we obtain $\|\psi_j^-\| = O(\exp(-c\epsilon^{-1/2}))$ for some $c > 0$.

Consider $f = c_1\psi_1^+ + c_2\psi_2^+$ for any $c \in \mathbf{R}^2$. Since $\phi_1 \perp \phi_2$, both have norm 1, and $\|\phi_j - \phi_j^+\| = O(\epsilon)$, it follows that $\|f\| = (1 + O(\epsilon))|c|$. Letting $H = -\partial_x^2 + (x - \epsilon^{-1})^2 - a(0) - b(z)$, by Lemma 3.2 we have $Hf = c_1\lambda_1\psi_1^+ + c_2\lambda_2\psi_2^+ + O(\epsilon^{1/2})|c|$, so $\|Hf\| \leq (4\theta + C\epsilon^{1/2})\|f\|$ for all $f$ in the two-dimensional space spanned by $\psi_1^+, \psi_2^+$. Thus by the minimax principle, $H$ has at least two eigenvalues in $[-1/2, 1/2]$ if $\theta, \epsilon$ are sufficiently small, a contradiction.

To prove existence of one small eigenvalue fix an eigenfunction $h$ of $-\partial_s^2 + s^2$ with eigenvalue $a(0)$. Setting $\phi_\epsilon(x) = h(x - \epsilon^{-1})$ and using the fact that $h$ is a Schwartz function, one obtains $\|B_{z,\epsilon}\phi_\epsilon\| = O(|(z,\epsilon)|)$ by expanding $q_\epsilon$ and its derivative about $x = \epsilon^{-1}$. The minimax principle then guarantees existence of an eigenvalue whose absolute value is $O(|(z,\epsilon)|)$. □

Fix a small constant $\theta > 0$ as in the preceding lemma. Substitute $x = y + \epsilon^{-1}$, so that

$$B_{z,\epsilon} = -\partial_y^2 + p_\epsilon^2(y) - [a(0) + b(z + \epsilon zy)]p_\epsilon'(y), \tag{3.4}$$

with

$$p_\epsilon(y) = (\epsilon k)^{-1}\left((1 + \epsilon y)^k - 1\right) = y + O(\epsilon y^2) + O(\epsilon^{k-1}y^k)$$
$$p_\epsilon'(y) = 1 + O(\epsilon|y| + \epsilon^{k-1}|y|^{k-1}).$$

Note that $p_\epsilon$ vanishes at $y = 0$ and at $y = -2\epsilon^{-1}$. For all small $(z, \epsilon)$ denote by $\lambda(z, \epsilon)$ the unique small (in absolute value) eigenvalue of $B_{z,\epsilon}$, and by $\phi = \phi(z, \epsilon)$ a corresponding eigenfunction of norm 1. The operator $(B_{z,\epsilon} - \lambda(z, \epsilon))^{-1}$ is well-defined as a bounded linear operator from the orthocomplement of $\phi(z, \epsilon)$ to $L^2(\mathbf{R})$.

**Lemma 3.5.** *Assume that $a(0) \in \{1, 3, 5, \ldots\}$. Then there exist $\delta > 0$, $C < \infty$ such that for all sufficiently small $(z, \epsilon)$ and all $0 < r \leq \delta$, for all $f \in L^2(\mathbf{R})$ orthogonal to $\phi(z, \epsilon)$,*

$$\int \left|(B_{z,\epsilon} - \lambda(z, \epsilon))^{-1}f(y)\right|^2 e^{r|y|}\, dy \leq C \int |f(y)|^2 e^{r|y|}\, dy.$$

*Proof.* Since $|\lambda(z, \epsilon)| \leq \theta$ and no other element of the spectrum of $B_{z,\epsilon}$ lies in $[-4\theta, 4\theta]$,

$$(B_{z,\epsilon} - \lambda(z, \epsilon))^{-1} = (2\pi i)^{-1}\int_{|\zeta - \lambda| = 2\theta}(B_{z,\epsilon} - \zeta)^{-1}\, d\zeta$$

as operators from the orthocomplement of $\phi(z, \epsilon)$ to $L^2$. Thus it suffices to establish the conclusion of the lemma for $(B_{z,\epsilon} - \zeta)^{-1}$ for all $\zeta$ on the contour of integration, uniformly in $\zeta$. The assumption that $f \perp \phi(z,\epsilon)$ is then no longer needed, as will be shown.

Let $f \in L^2$ be given and set $g = (B_{z,\epsilon} - \zeta)^{-1} f \in L^2$. Repeating the reasoning in the proof of Lemma 3.2 and exploiting the assumption that $a(0) > 0$ and hence $(p_\epsilon^2(y) - [a(0) + b(z+\epsilon zy)]p'_\epsilon(y))$ has a strictly positive lower bound for $|y + 2\epsilon^{-1}| \leq C_2$, uniformly for $(z,\epsilon)$ sufficiently close to 0, one obtains for each $N < \infty$

$$\int |g|^2 e^{\min(r|y|,N)}\, dy$$
$$\leq C \int |(B_{z,\epsilon} - \zeta)g|^2\, e^{\min(r|y|,N)}\, dy + C \int_{|y| \leq C_2} |g|^2\, e^{\min(r|y|,N)}\, dy$$
$$\leq C \int |f|^2\, e^{\min(r|y|,N)}\, dy + C \int |g|^2\, dy$$
$$\leq C \int |f|^2\, e^{\min(r|y|,N)}\, dy + C\theta^{-1} \int |f|^2\, dy$$
$$\leq C \int |f|^2\, e^{r|y|}\, dy,$$

uniformly in $z, \epsilon, \zeta, N$ provided $z, \epsilon$ are small and $|\zeta - \lambda| = 2\theta$. Taking the limit as $N \to \infty$ concludes the proof. $\square$

The same analysis and conclusions hold for $\epsilon < 0$, as well, provided that $|\epsilon|$ is sufficiently small.

3.2. **Perturbation expansions and smooth dependence of eigenvalues.** If $b$ does not vanish to infinite order at 0 then the operators $B_{z,\epsilon}$ do not depend smoothly on $\epsilon$ *uniformly* as $\epsilon \to 0$. More precisely, the norm of the formal derivative $\partial^n[b(z+\epsilon zy)p'_\epsilon(y)]/\partial \epsilon^n$, as an operator from the domain of $B_{z,\epsilon}$ to $L^2$, tends to infinity like some negative power of $|\epsilon|$ once $n$ is sufficiently large. Nevertheless $\lambda(z,\epsilon)$ will be shown to extend to a $C^\infty$ function in a neighborhood of $0 \in \mathbf{R}^2$. The next lemma is one ingredient in the proof. Denote by $C^{N,1}$ the class of $N$ times continuously differentiable functions whose partial derivatives of order $N$ are all Lipschitz continuous.

**Lemma 3.6.** *Let $\Omega \subset \mathbf{R}^n$ be an open set, $F : \Omega \mapsto \mathbf{C}$ a function, and $N \geq 0$ an integer. Suppose there exists $C' < \infty$ such that for each $x \in \Omega$ there exists a polynomial $P_x$ of degree not exceeding $N$ such that for all $y \in \Omega$,*

$$|F(y) - P_x(y)| \leq C'|y - x|^{N+1},$$

and all coefficients of each $P_x$ are bounded in modulus by $C'$. Then for any relatively compact open $\Omega' \subset \Omega$, $F$ belongs to $C^{N,1}(\Omega')$ with norm bounded by a constant depending only on $C', N$ and the distance from $\Omega'$ to the complement of $\Omega$. Moreover at each point $x$, $P_x$ is the Taylor polynomial of degree $N$ for $F$ at $x$.

*Proof.* Assume $|x - x'|$ is at most one third of the distance from $x$ to the boundary of $\Omega$. Then $|P_x(y) - P_{x'}(y)| \leq C|x - x'|^{N+1}$ whenever $|y - x| \leq 2|x - x'|$. For each $N$ there exists a constant $A_N < \infty$ such that for any polynomial $Q$ of degree at most $N$, for any $|\alpha|$,
$$|\partial^\alpha Q(0)| \leq A_n \sup_{|w| \leq 1} |Q(w)|.$$
Applying this to $Q(w) = P_x(x + |x - x'|w) - P_{x'}(x + |x - x'|w)$ yields $|\partial_y^\alpha P_x(y) - \partial_y^\alpha P_{x'}(y)| \leq C|x - x'|^{N+1-|\alpha|}$ for all $0 \leq |\alpha| \leq N$. Define
$$F_\alpha(x) = \partial_y^\alpha P_x(y)\Big|_{y=x}.$$
Then each $F_\alpha$ is Lipschitz continuous, for
$$|F_\alpha(x) - F_\alpha(x')| \leq |\partial_y^\alpha P_x(x) - \partial_y^\alpha P_{x'}(x)| + |\partial_y^\alpha P_{x'}(x) - \partial_y^\alpha P_{x'}(x')|$$
$$\leq C|x - x'|.$$

It follows that on any compact subset $K$ of $\Omega$, each $F_\alpha$ is bounded by a constant depending only on $C', N, \Omega, K$. Setting $P_x^{(M)}(y) = \sum_{|\alpha| \leq M} F_\alpha(x)(y - x)^\alpha/\alpha!$ for any $0 \leq M \leq N$, we find that the hypotheses of the Lemma with $N$ replaced by $M$ are also satisfied by the polynomials $P_x^{(M)}$. It then follows by induction on $M$ that $F \in C^{M,1}$. $\square$

The formula (3.4) for $B_{z,\epsilon}$ makes sense for $\epsilon < 0$ as well as for $\epsilon > 0$, and by continuity extends to $\epsilon = 0$ in such a way that as a map from the Schwartz class to $L^2$, $B_{z,\epsilon}$ depends smoothly on $z, \epsilon$. The above analysis applies also for $\epsilon < 0$ and demonstrates existence of a unique small eigenvalue $\lambda(z, \epsilon)$. Denote by $\phi(z, \epsilon)$ an associated eigenfunction of norm 1. By the proof of Lemmas 3.2 and 3.5, there exist $r, C \in \mathbf{R}^+$ such that
$$\int \phi(z, \epsilon)^2(y)\, e^{r|y|}\, dy \leq C \tag{3.5}$$
for all $z, \epsilon$ in a neighborhood of 0.

**Lemma 3.7.** *There exist $\delta > 0$ and bounded coefficients $\Lambda_j(\zeta, \epsilon)$ such that for every positive integer $N$, for every $z, \zeta, \epsilon \in [-\delta, \delta]$,*
$$|\lambda(z, \epsilon) - \sum_{j=0}^N \Lambda_j(\zeta, \epsilon)(z - \zeta)^j| \leq C_N |z - \zeta|^{N+1}.$$

*Proof.* Write Taylor expansions
$$B_{z,\epsilon} \sim \sum_{j=0}^{\infty} \beta_j(\zeta,\epsilon)(z-\zeta)^j,$$
where $\beta_0(\zeta,\epsilon) = B_{\zeta,\epsilon}$ and for $j \geq 1$,
$$\beta_j(\zeta,\epsilon) = -(1+\epsilon y)^j p'_\epsilon(y) b^{(j)}(\zeta + \epsilon\zeta y)/j!;$$
$\beta_j$ denotes both a function of $y$ and the operator defined by multiplication by that function. Fix $N \geq 0$ and write
$$\Lambda(z,\epsilon) = \sum_{j=0}^{N} \Lambda_j(\zeta,\epsilon)(z-\zeta)^j,$$
$$\psi(z,\epsilon) = \sum_{j=0}^{N} \psi_j(\zeta,\epsilon)(z-\zeta)^j$$
with $\Lambda_0(\zeta,\epsilon) = \lambda(\zeta,\epsilon)$, $\psi_0(\zeta,\epsilon) = \phi(\zeta,\epsilon)$, where $\Lambda_j$ and $\psi_j$ are to be determined for $j \geq 1$ by solving the equation
$$B_{z,\epsilon}\psi(z,\epsilon) = \Lambda(z,\epsilon)\psi(z,\epsilon) + O(|z-\zeta|^{N+1}). \tag{3.6}$$

Equating like powers of $z - \zeta$ in this equation yields
$$[B_{\zeta,\epsilon} - \lambda(\zeta,\epsilon)]\psi_n(\zeta,\epsilon) = -\sum_{j=1}^{n}[\beta_j(\zeta,\epsilon) - \Lambda_j(\zeta,\epsilon)]\psi_{n-j}(\zeta,\epsilon) \tag{3.7}$$
for all $1 \leq n \leq N$. The unknowns $\Lambda_n, \psi_n$ are determined by induction on $n$; if the right-hand side is given and belongs to $L^2(\mathbf{R})$ then a necessary and sufficient condition for the existence of a solution $\psi_n \in L^2(\mathbf{R})$ is that
$$0 = \langle \phi(\zeta,\epsilon), \sum_{j=1}^{n}[\beta_j(\zeta,\epsilon) - \Lambda_j(\zeta,\epsilon)]\psi_{n-j}(\zeta,\epsilon)\rangle,$$
which, since $\langle \phi(\zeta,\epsilon), \phi(\zeta,\epsilon)\rangle = 1$, may be rewritten as
$$\Lambda_n(\zeta,\epsilon) = \langle \beta_n(\zeta,\epsilon)\phi(\zeta,\epsilon), \phi(\zeta,\epsilon)\rangle$$
$$+ \sum_{j=1}^{n-1} \langle \phi(\zeta,\epsilon), [\beta_j(\zeta,\epsilon) - \Lambda_j(\zeta,\epsilon)]\psi_{n-j}(\zeta,\epsilon).\rangle. \tag{3.8}$$

For $n = 1$ this last sum is vacuous, and the equation reads
$$\Lambda_1(\zeta,\epsilon) = \langle \beta_1(\zeta,\epsilon)\phi(\zeta,\epsilon), \phi(\zeta,\epsilon)\rangle.$$

Once $\Lambda_j, \psi_j$ are defined for all $0 \le j < n$, (3.8) determines $\Lambda_n(\zeta, \epsilon)$ uniquely in terms of those $\Lambda_j, \psi_j$. (3.7) then uniquely determines $\psi_n \in L^2$, provided that the right-hand side in (3.7) does belong to $L^2$. This last point requires some justification, however, since the $\beta_j$ are not bounded operators on $L^2$.

Fix an infinite sequence of small exponents $r_0 > r_1 > \cdots > 0$, all satisfying the conclusions of Lemma 3.5. (3.5) guarantees in particular that $\psi_0(\zeta, \epsilon) = \phi(\zeta, \epsilon)$ is bounded in $L^2(\mathbf{R}, \exp(r_0|y|dy))$, uniformly for $(\zeta, \epsilon)$ near 0. $\beta_1(\zeta, \epsilon)$ is multiplication by a function bounded by $C(1 + |y|)^M$ for some finite $M$, uniformly in $(\zeta, \epsilon)$, so the right-hand side of (3.7) belongs to $L^2(\mathbf{R}, \exp(r_1|y|dy))$, still uniformly in $(\zeta, \epsilon)$. By induction on $n$ and by Lemma 3.5, the unique solution $\psi_n \in L^2(\mathbf{R})$ of (3.7) belongs to $L^2(\mathbf{R}, \exp(r_n|y|dy))$, uniformly in $(\zeta, \epsilon)$. Boundedness of the coefficients $\Lambda_j(\zeta, \epsilon)$ follows in the same way.

It remains to verify that

$$\lambda(z, \epsilon) = \sum_{j=0}^{N} \Lambda_j(\zeta, \epsilon)(z - \zeta)^j + O(|z - \zeta|^{N+1}).$$

Setting $\psi(z, \epsilon) = \sum_{j=0}^{N} \psi_j(\zeta, \epsilon)(z - \zeta)^j$ and $\Lambda(z, \epsilon) = \sum_{j=0}^{N} \Lambda_j(\zeta, \epsilon)(z - \zeta)^j$, we have

$$(B_{z,\epsilon} - \Lambda(z, \epsilon))\psi(z, \epsilon) = O(|z - \zeta|^{N+1})$$

in $L^2$ norm, by construction, and $\|\psi(z, \epsilon)\| = 1 + O(|z - \zeta|) \ge 1/2$ provided that $|z - \zeta|$ is sufficiently small. Since $B_{z,\epsilon}$ is selfadjoint, this forces the distance from $\Lambda(z, \epsilon)$ to the spectrum of $B_{z,\epsilon}$ to be $O(|z - \zeta|^{N+1})$. But $\Lambda_0(\zeta, \epsilon) = \lambda(\zeta, \epsilon)$ by definition and the latter is small, so $|\Lambda(z, \epsilon)| \le 2\theta$ for all $(z, \epsilon)$ sufficiently close to 0. Since $B_{z,\epsilon}$ has discrete spectrum and $\lambda(z, \epsilon)$ is its only eigenvalue in $[-4\theta, 4\theta]$, this forces $|\Lambda(z, \epsilon) - \lambda(z, \epsilon)| = O(|z - \zeta|^{N+1})$. $\square$

$B_{z,\epsilon}$, in the form of equation (3.4), extends to $\epsilon = 0$ as a $C^\infty$ function of all $(z, \epsilon)$ in a neighborhood of $0 \in \mathbf{R}^2$. The same reasoning as in the proof of Lemma 3.7 therefore yields bounded coefficients $\Lambda_\alpha(\zeta, \epsilon)$ satisfying

$$\lambda(z, e) = \sum_{0 \le |\alpha| \le N} \Lambda_\alpha(\zeta, \epsilon)((z - \zeta), (e - \epsilon))^\alpha + O(|(z, e) - (\zeta, \epsilon)|^{N+1})$$

From Lemma 3.6 we then draw the following conclusion.

**Corollary 3.8.** *For $(z, \epsilon)$ in a sufficiently small neighborhood of the origin, the unique small eigenvalue $\lambda(z, \epsilon)$ of $B_{z,\epsilon}$ is a $C^\infty$ function of $(z, \epsilon)$.*

**Corollary 3.9.** *As a function of $z$, in some neighborhood of $0 \in \mathbf{R}^2$, $\lambda(z, \epsilon) = -b(z) + O(\epsilon)$ in the $C^N$ norm for any $N$.*

*Proof.* The set of all eigenvalues of $B_{z,0} = -\partial_y^2 + y^2 - a(0) - b(z)$ is the set of all numbers $\lambda - [a(0) + b(z)]$ such that $\lambda \in \{1, 3, 5, \ldots\}$. Since $a(0) \in \{1, 3, 5, \ldots\}$, $-b(z)$ is therefore the unique small eigenvalue when $\epsilon = 0$. The result for small $\epsilon \neq 0$ then follows from the preceding corollary. □

**Corollary 3.10.** *Suppose that $a(0) \in \{1, 3, 5, \ldots\}$ and that $a^{(m)}(0) \neq 0$ for some $m \geq 1$. Then*
$$\frac{\partial^m \lambda(z, \epsilon)}{\partial z^m} \neq 0 \text{ for all } (z, \epsilon) \text{ sufficiently close to } 0.$$

All this reasoning applies equally well when $\tau < 0$ and/or $a(0)$ belongs to the set of all negative odd integers. It applies also for odd $k$ with a simplification, since $q_\epsilon$ then has only a single zero. For $k$ even, $q_\epsilon$ has no zeros when $\eta\tau > 0$, resulting in the strong bound $\|A_{\eta,\tau}f\| \geq c|\tau|^{2/(k+1)}\|f\|$ of Lemma 2.2. When $k$ is odd there is no distinction between the cases $\eta\tau > 0$ and $\eta\tau < 0$; $q_\epsilon$ has one zero in each case. When $k = 1$ there are additional simplifications, since $q_\epsilon$ is then a linear function of $y$, but the same reasoning still applies.

3.3. **Finite measure of exceptional parameter sets.** Recall that $\{\mu_j\}$ denote the eigenvalues of $A_{\eta,\tau}$. The next result is Lemma 3.4 of [2], where a proof may be found. Denote by $B^n$ the closed unit ball in $\mathbf{R}^n$.

**Lemma 3.11.** *Suppose that $n, m \geq 1$, that $f \in C^{m+1}(B^n)$, and that there exists a multi-index $\alpha$ satisfying $0 \leq |\alpha| \leq m$ such that for every $y \in B^n$, $\partial^\alpha f/\partial x^\alpha(y) \neq 0$. Then there exists a constant $C < \infty$ such that for every $\delta > 0$,*
$$|\{y \in B^n : |f(y)| \leq \delta\}| \leq C\delta^{1/m}. \tag{3.9}$$

**Lemma 3.12.** *For any $k \geq 1$, if $a(0) \in \{\pm 1, \pm 3, \pm 5, \ldots\}$ and $a^{(m)}(0) \neq 0$ for some $m \geq 1$, then there exists $M < \infty$ such that*
$$\left|\{(\eta, \tau) \in \mathcal{C}_k \min_j |\mu_j(\eta, \tau)| \leq |\tau|^{-M}\}\right| < \infty.$$

*Proof.* Consider first the case where $k$ is even. For $\eta\tau > 0$ one has $\|A_{\eta,\tau}f\| \geq c|\tau|^{2/(k+1)}\|f\|$ for all $f, \eta, \tau$ by Lemma 2.2, so the exceptional set in question is empty for any $M > 1$. The case $\tau < 0$ reduces to $\tau > 0$ by the change of variables $(x, y, t) \mapsto (x, y, -t)$, so we may assume the latter.

The Jacobian determinant for the change of variables $(\eta, \tau) \mapsto (z, \epsilon)$ introduced above is
$$\left|\frac{\partial(\eta, \tau)}{\partial(z, \epsilon)}\right| = c|\epsilon^{-5} z^{-k-3}| \leq C\tau^R$$
for some $C, R \in \mathbf{R}^+$. Consider $S_q = \{(\eta, \tau) \in \mathcal{C}_k : 2^q \leq \tau \leq 2^{q+1}\}$ where $q$ is an arbitrary large positive integer. $S_q$ is mapped into an arbitrarily small neighborhood of 0 as $q \to \infty$.

Since $A_{\eta,\tau}$ is unitarily equivalent to $|\epsilon z|^{-2} B_{z,\epsilon}$, $\min_j |\mu_j(\eta, \tau)|$ equals $|\epsilon z|^{-2}|\lambda(z, \epsilon)|$. The quantity $|\epsilon z|^{-2}$ equals a constant times $|\eta|^{(k-1)/k}|\tau|^{1/k}$, and both $|\eta|$ and $|\tau|$ are bounded below by a positive constant when $(\eta, \tau) \in \mathcal{C}_k$. This constant may be taken to be at least 1, by choosing $\tau_0$ to be sufficiently large. Thus $\min_j |\mu_j(\eta, \tau)| \geq |\lambda(z, \epsilon)|$ for all $(\eta, \tau) \in \mathcal{C}_k$.

Let $n \geq 1$ be an index for which $a^{(n)}(0) \neq 0$. Fix an exponent $M > nR$. By Corollary 3.10 there exists $\delta > 0$ such that $\partial^n \lambda(z, \epsilon)/\partial z^n \neq 0$ for all $|(z, \epsilon)| \leq 2\delta$, so for each $|\epsilon| < \delta$,
$$\left|\{|z| < \delta : |\lambda(z, \epsilon)| \leq 2^{-qM}\}\right| \leq C 2^{-qM/n}$$
by Lemma 3.11. One has $\epsilon^2 \leq k^{-1} \gamma_0^{-(k+1)/k}$, so this will be satisfied for all $(\eta, \tau) \in \mathcal{C}_k$ provided that $\gamma_0$ is chosen to be sufficiently large. Therefore
$$\left|\{(\eta, \tau) \in \mathcal{C}_k : 2^q \leq \tau \leq 2^{q+1} \text{ and } \min_j |\mu_j(\eta, \tau)| \leq |\tau|^{-M}\}\right| \leq C 2^{qR} 2^{-qM/n}.$$
Summing over $q$ yields the desired conclusion.

The reasoning for odd $k$ is the same. □

## 4. The case $a(0) \in \{\pm 3, \pm 5, \dots\}$ for $k > 1$

We continue to assume that $(\eta, \tau) \in \mathcal{C}_k$. A different analysis is required in this case if $a^{(m)}(0) = 0$ for all $m \geq 1$, for the case $a(0) = \pm 1$ must be distinguished from the other exceptional cases. Define new variables $(z, \epsilon)$ in terms of $(\eta, \tau)$ as above, let $B_{z,\epsilon}$ be the operator defined in (3.4) and for small $(z, \epsilon)$ let $\lambda(z, \epsilon)$ continue to denote its unique small eigenvalue. An asymptotic expansion
$$\lambda(z, \epsilon) \sim \sum_{j \geq 0} \Lambda_j(z) \epsilon^j$$
has already been established, with $\Lambda_j \in C^\infty$ in a neighborhood of the origin and $\Lambda_0(z) = -b(z)$. Set $\Lambda_j = \Lambda_j(0)$, so that $\lambda(0, \epsilon) \sim \sum \Lambda_j \epsilon^j$.

**Lemma 4.1.** *Assume that $a(0) = 2n + 1$ for some integer $n \geq 0$. Then*
$$\Lambda_2 = (k-1)n(n+1)/2.$$

The proof will show that $\Lambda_1 = 0$, so it is necessary to pass to the second coefficient in the expansion. Although it is easy to see without elaborate calculation that every perturbation coefficient $\Lambda_j$ must vanish when either $k = 1$ or $n = 0$, we can offer no simple or conceptual explanation for the nonvanishing of $\Lambda_2$ when $(k-1)n \neq 0$. Before presenting the calculations we record their consequence.

**Corollary 4.2.** *If $k > 1$ and $a(0) \in \{3, 5, \dots\}$ then $\partial^2 \lambda(z, \epsilon)/\partial \epsilon^2 \neq 0$ in some neighborhood of $0 \in \mathbf{R}^2$.*

With Lemma 4.1 in hand, the proof of the next lemma is parallel to that of Lemma 3.12 and is therefore omitted.

**Lemma 4.3.** *If $k > 1$ and $a(0) \in \{\pm 3, \pm 5, \dots\}$ then there exists $M < \infty$ such that*
$$\left|\{(\eta, \tau) \in \mathcal{C}_k : \min_j |\mu_j(\eta, \tau)| \leq |\tau|^{-M}\}\right| < \infty.$$

**Proof of Lemma 4.1.** Begin with the Taylor expansion
$$\begin{aligned} p_\epsilon(y) &= (\epsilon k)^{-1}((1 + \epsilon y)^k - 1) \\ &= y + [(k-1)y^2/2]\epsilon + [(k-1)(k-2)y^3/6]\epsilon^2 + O(\epsilon^3), \\ p'_\epsilon(y) &= 1 + [(k-1)y]\epsilon + [(k-1)(k-2)y^2/2]\epsilon^2 + O(\epsilon^3). \end{aligned}$$

Thus
$$\begin{aligned} &p_\epsilon^2(y) - a(0)p'_\epsilon(y) \\ &= y^2 - a(0) + \left[(k-1)y^3 - a(0)(k-1)y\right]\epsilon \\ &\quad + \left[\left(\frac{k-1}{2}\right)^2 y^4 + \frac{(k-1)(k-2)}{3}y^4 - a(0)\frac{(k-1)(k-2)}{2}y^2\right]\epsilon^2 + O(\epsilon^3). \end{aligned}$$

Expanding $B_{0,\epsilon} \sim \sum_j \beta_j \epsilon^j$, one has
$$\begin{aligned} \beta_0 &= H_n = -\partial_y^2 + y^2 - (2n+1) \\ \beta_1 &= (k-1)(y^3 - (2n+1)y) \\ \beta_2 &= (k-1)\left[\frac{7k-11}{12}y^4 - (2n+1)\frac{k-2}{2}y^2\right]. \end{aligned}$$

By the same reasoning as in section 3 there exist $\psi_i \in L^2(\mathbf{R}, \exp(r|y|)dy)$ for some $r > 0$ and scalars $\Lambda_i$ such that $H_n \psi_0 = 0$ and
$$B_{0,\epsilon}(\sum_{i=0}^{2} \psi_i) = (\sum_{i=0}^{2} \Lambda_i \epsilon^i)(\sum_{i=0}^{2} \psi_i) + O(\epsilon^3),$$

where $\psi_0$ does not vanish identically. Consequently, as in section 3, we have

$$\lambda(0, \epsilon) = \sum_{i=0}^{2} \Lambda_i \epsilon^i + O(\epsilon^3).$$

By (3.8),

$$\Lambda_1 = \langle \beta_1 \psi_0, \psi_0 \rangle = (k-1) \int_{\mathbf{R}} (y^3 - (2n+1)y)\psi_0^2(y)\, dy = 0$$

for any $n$, because the Hermite eigenfunction $\psi_0$ is either even or odd, hence its square is even. By (3.7),

$$\psi_1 = -H_n^{-1}(\beta_1 \psi_0) = -(k-1)H_n^{-1}((y^3 - (2n+1)y)\psi_0).$$

Then (3.8) gives

$$\Lambda_2 = \langle \psi_0, \beta_2 \psi_0 + \beta_1 \psi_1 \rangle$$
$$= \left\langle \psi_0,\ (k-1)\left(\frac{7k-11}{12}y^4 - (2n+1)\frac{k-2}{2}y^2\right)\psi_0 \right\rangle$$
$$- \left\langle \psi_0,\ (k-1)(y^3 - (2n+1)y)H_n^{-1}(k-1)(y^3 - (2n+1)y)\psi_0 \right\rangle$$

so that

$$(k-1)^{-1}\Lambda_2 = \frac{7k-11}{12}\|y^2\psi_0\|^2 - (2n+1)\frac{k-2}{2}\|y\psi_0\|^2$$
$$- (k-1)\left\langle (y^3 - (2n+1)y)\psi_0,\ H_n^{-1}\left((y^3 - (2n+1)y)\psi_0\right) \right\rangle. \tag{4.1}$$

We next recall certain formulas concerning Hermite eigenfunctions and their derivations. Let

$$H = -\partial_y^2 + y^2$$
$$h_0(y) = c_0 e^{-y^2/2}$$

with $c_0$ chosen so that $\|h_0\| = 1$. Then $Hh_0 = h_0$. Inductively define

$$h_{q+1} = [2(q+1)]^{-1/2}(-\partial_y + y)h_q.$$

Then $h_q$ is an eigenfunction of $H$ with eigenvalue $2q+1$. Moreover $\|h_q\| = 1$ for all $q$ because

$$\|h_{q+1}\|^2 = [2(q+1)]^{-1}\langle (\partial_y + y)(-\partial_y + y)h_q,\ h_q \rangle$$
$$= [2(q+1)]^{-1}\langle (H+1)h_q,\ h_q \rangle$$
$$= [2(q+1)]^{-1}(2q+2)\|h_q\|^2.$$

The lowering identity is

$$(\partial_y + y)h_q = (2q)^{1/2}h_{q-1}, \tag{4.2}$$

because $-\partial_y + y$ is injective and $(-\partial_y + y)(\partial_y + y)h_q = (H-1)h_q = 2qh_q$ while $(-\partial_y + y)(2q)^{1/2}h_{q-1} = (2q)^{1/2}(2q)^{1/2}h_q$ by definition of $h_q$. Combining the lowering identity with the definition of $h_{q+1}$ gives

$$\begin{aligned} yh_q &= 2^{-1}(\partial_y + y)h_q + 2^{-1}(-\partial_y + y)h_q \\ &= 2^{-1/2}q^{1/2}h_{q-1} + 2^{-1/2}(q+1)^{1/2}h_{q+1}. \end{aligned}$$

Iterating this last formula gives

$$y^2 h_q = 2^{-1}(q+1)^{1/2}(q+2)^{1/2}h_{q+2} + 2^{-1}(2q+1)h_q + 2^{-1}q^{1/2}(q-1)^{1/2}h_{q-2}.$$

Iterating once more yields

$$\begin{aligned} y^3 h_q &= 2^{-3/2}\left[(q+3)^{1/2}(q+2)^{1/2}(q+1)^{1/2}\right]h_{q+3} \\ &\quad + 2^{-3/2}\left[(q+2)(q+1)^{1/2} + (2q+1)(q+1)^{1/2}\right]h_{q+1} \\ &\quad + 2^{-3/2}\left[(2q+1)q^{1/2} + q^{1/2}(q-1)\right]h_{q-1} \\ &\quad + 2^{-3/2}\left[q^{1/2}(q-1)^{1/2}(q-2)^{1/2}\right]h_{q-3} \\ &= 2^{-3/2}\left[(q+3)(q+2)(q+1)\right]^{1/2}h_{q+3} \\ &\quad + 2^{-3/2}3(q+1)^{3/2}h_{q+1} \\ &\quad + 2^{-3/2}3q^{3/2}h_{q-1} \\ &\quad + 2^{-3/2}\left[q(q-1)(q-2)\right]^{1/2}h_{q-3}. \end{aligned}$$

Therefore

$$\|yh_q\|^2 = 2^{-1}q + 2^{-1}(q+1) = (2q+1)/2,$$

and

$$\begin{aligned} \|y^2 h_q\|^2 &= 2^{-2}\left[(q+1)(q+2) + (2q+1)^2 + q(q-1)\right] \\ &= 4^{-1}\left[q^2 + 3q + 2 + 4q^2 + 4q + 1 + q^2 - q\right] \\ &= \frac{3}{4}(2q^2 + 2q + 1). \end{aligned}$$

Also
$$(y^3 - (2q+1)y)h_q = 2^{-3/2}\left[(q+3)(q+2)(q+1)\right]^{1/2} h_{q+3}$$
$$+ 2^{-3/2}\left[3(q+1)^{3/2} - 2(2q+1)(q+1)^{1/2}\right] h_{q+1}$$
$$+ 2^{-3/2}\left[3q^{3/2} - 2(2q+1)q^{1/2}\right] h_{q-1}$$
$$+ 2^{-3/2}\left[q(q-1)(q-2)\right]^{1/2} h_{q-3}.$$

$H_q$ was defined to be $H - (2q+1)$, so $H_q^{-1}h_p = 2^{-1}(p-q)^{-1}h_p$ for all $p \neq q$. Therefore

$$3 \cdot 2^4 \left\langle (y^3 - (2q+1)y)h_q, \; H_q^{-1}(y^3 - (2q+1)y)h_q \right\rangle$$
$$= [(q+3)(q+2)(q+1)] + 3(q+1)[3(q+1) - 2(2q+1)]^2$$
$$\quad - 3q[3q - 2(2q+1)]^2 - [q(q-1)(q-2)]$$
$$= (q^3 + 6q^2 + 11q + 6) + 3(q+1)(-q+1)^2$$
$$\quad - 3q(-q-2)^2 - (q^3 - 3q^2 + 2q)$$
$$= (q^3 + 6q^2 + 11q + 6) + (3q^2 - 3q^2 - 3q + 3)$$
$$\quad + (-3q^3 - 12q^2 - 12q) + (-q^3 + 3q^2 - 2q)$$
$$= -6q^2 - 6q + 9.$$

Combining all these ingredients yields a formula for $\Lambda_2$.

$$16(k-1)^{-1}\Lambda_2$$
$$= 16\frac{7k-11}{12}\frac{3}{4}(2n^2 + 2n + 1)$$
$$\quad - 16(2n+1)\frac{k-2}{2}\frac{2n+1}{2} - 16\frac{k-1}{48}(-6n^2 - 6n + 9)$$
$$= (7k-11)(2n^2 + 2n + 1) - 4(k-2)(4n^2 + 4n + 1) + (k-1)(2n^2 + 2n - 3)$$
$$= (n^2 + n)(14k - 22 - 16k + 32 + 2k - 2) + (7k - 11 - 4k + 8 - 3k + 3)$$
$$= 8n(n+1). \quad \square$$

5. The nonperturbative parameter regime $\mathcal{D}_k$

For $(\eta, \tau) \in \mathcal{D}_k$ define
$$\epsilon = |\tau|^{-1/(k+1)}, \qquad w = \operatorname{sgn}(\tau)\eta|\tau|^{-1/(k+1)}$$
and
$$D_{w,\epsilon} = -\partial_x^2 + (x^k + w)^2 - \operatorname{sgn}(\tau)[a(0) + b(\epsilon x)]kx^{k-1}.$$

Then $A_{\eta,\tau}$ is unitarily equivalent to $|\tau|^{2/(k+1)} D_{w,\epsilon}$, via the substitution $s = \epsilon x$. $D_{w,\epsilon}$ is essentially selfadjoint for each $k \geq 1$ and each $(w, \epsilon) \in \mathbf{R}^2$, has compact resolvent, and its spectrum consists of a sequence of real eigenvalues $\lambda_0(w, \epsilon) < \lambda_1(w, \epsilon) < \ldots$ tending to $+\infty$. All eigenspaces are one dimensional.

$(\eta, \tau) \in \mathcal{D}_k$ if and only if $0 < \epsilon \leq \tau_0^{-1/(k+1)} \ll 1$ and $|w| \leq \gamma_0 < \infty$. The analysis for $\mathcal{D}_k$ differs from that for $\mathcal{C}_k$ in that $\mathcal{D}_k$ is not a perturbative regime; we are not able to analyze $D_{w,\epsilon}$ by showing that it is close to a better understood operator. In particular, although the definition of $\mathcal{D}_k$ requires $\epsilon$ to be close to 0, the constant $\gamma_0$ must be taken to be sufficiently large in order for the analysis of $\mathcal{C}_k$ above to succeed. Thus $\mathcal{D}_k$ must encompass the case where $\epsilon = 0$ but $w$ is bounded by a large constant. We will instead derive information for bounded $w$ by analytic continuation from the case of large $w$, which has already been treated by perturbative techniques.

**Lemma 5.1.** *The eigenvalues $\lambda_n(w, \epsilon)$ are $C^\infty$ functions of $(w, \epsilon) \in \mathbf{R}^2$ and are real analytic functions of $w$, uniformly for all $\epsilon$ in any compact subset of $\mathbf{R}$.*

*Proof.* Formally $D_{w,\epsilon}$ depends holomorphically on $w \in \mathbf{C}$, for each fixed $\epsilon$. We claim that $D_{w,\epsilon}$ is a bounded operator from the domain of $D_{0,0}$ to $L^2$ and satisfies
$$\|D_{w,\epsilon} f\| \leq C\|D_{0,0} f\| + C\|f\|$$
for all $f$ in the domain of $D_{0,0}$, uniformly for $(w, \epsilon)$ in any compact subset of $\mathbf{C} \times \mathbf{R}$. Consequently $w \mapsto D_{w,\epsilon}$ is an analytic family of operators in the sense of Kato [8], and since the spectrum consists entirely of eigenvalues associated to one dimensional eigenspaces, the theory of such families guarantees holomorphic dependence of the eigenvalues on $w$ and their extension to entire holomorphic functions of $w \in \mathbf{C}$, given that no two ever coalesce, which we already know to be true.

To prove the inequality, it suffices to consider any $f \in C_0^2$. Then
$$\langle D_{0,0} f, f \rangle = \|\partial_x f\|^2 + \|x^k f\|^2 \pm a(0) k \int x^{k-1} |f|^2 \, dx,$$
so
$$\|\partial_x f\| + \|x^k f\| \leq C\|D_{0,0} f\| + C\|f\|.$$
Consequently
$$\|(D_{w,\epsilon} - D_{0,0}) f\| = \left\| \left( 2w x^k + w^2 \pm k b(\epsilon x) x^{k-1} \right) f \right\|$$
$$\leq C_{w,\epsilon} (\|D_{0,0} f\| + \|f\|),$$
as desired. The same reasoning yields an inequality
$$\|(D_{w,\epsilon} - D_{w',\epsilon'}) f\| \leq C|(w, \epsilon) - (w', \epsilon')| [\|D_{w,\epsilon} f\| + \|f\|] \tag{5.1}$$
provided that $w, \epsilon, w', \epsilon'$ are assumed to lie in any fixed bounded region.

The method of proof in section 3 establishes that each $\lambda_n$ is a $C^\infty$ function of $(w, \epsilon)$, once $\lambda_n$ is known to be a continuous function. Each $D_{w,\epsilon}$ has discrete spectrum consisting of eigenvalues associated to one dimensional eigenspaces, and has the same domain as $D_{0,0}$. Given any $(w, \epsilon)$ and any compact set $K$ disjoint from the spectrum of $D_{w,\epsilon}$, (5.1) guarantees that $K$ is also disjoint from the spectrum of $D_{w',\epsilon'}$ for all $(w', \epsilon')$ sufficiently close to $(w, \epsilon)$. On the other hand, given any eigenvalue $\lambda$ of $D_{w,\epsilon}$, fix a circle $\Gamma$ centered at $\lambda$ such that all other eigenvalues of $D_{w,\epsilon}$ lie in the exterior of $\Gamma$. Then $P_{w',\epsilon'} = (2\pi i)^{-1} \oint_\Gamma (D_{w',\epsilon'} - z)^{-1} dz$ is a projection onto the direct sum of all eigenspaces of $D_{w',\epsilon'}$ associated to eigenvalues belonging to the disk bounded by $\Gamma$. By (5.1), $P_{w,\epsilon} - P_{w',\epsilon'} = O(|(w, \epsilon) - (w', \epsilon')|)$. Thus $P_{w',\epsilon'}$ must have rank one for all $(w', \epsilon')$ sufficiently close to $(w, \epsilon)$, so $D_{w',\epsilon'}$ has a unique eigenvalue in the interior of $\Gamma$. Taking $\Gamma$ to have arbitrarily small radius completes the proof. □

**Lemma 5.2.** *Let $k > 1$. If $k$ is even, or if $a(0) \notin \{\pm 1\}$, or if $\mathrm{sgn}\,(\tau)a(0) = -1$, then for each index $n$ the function $\mathbf{C} \ni w \mapsto \lambda_n(w, 0)$ does not vanish identically.*

*Proof.* Suppose first that $k$ is even, and consider the case where $w \in \mathbf{R}$ is positive and large. For any $f \in C_0^2$,

$$\langle D_{w,0} f, f \rangle \geq \|\partial_x f\|^2 + \|(x^k + w)f\|^2 - C \int |x|^{k-1} |f(x)|^2 \, dx$$
$$\geq \int \left( x^{2k} + w^2 - C|x|^{k-1} \right) |f|^2$$
$$\geq \frac{1}{2} w^2 \|f\|^2$$

for large $w$. Thus for every $n$, $\lambda_n(w, 0) \to +\infty$ as $w \to +\infty$.

Suppose next that $a(0) \notin \{\pm 1\}$ and $k$ is odd, and consider the case of large negative $w$. Set $\sigma = (k-1)/2k < 1/2$, and substitute $x = |w|^{-\sigma} y$, $|w| = \delta^{-1/(1-\sigma)}$ to obtain

$$D_{w,\epsilon} = |w|^{2\sigma} \left( -\partial_y^2 + p_\delta^2(y) - \mathrm{sgn}\,(\tau)[a(0) + b(\epsilon y)] p_\delta'(y) \right),$$

with $p_\delta(y) = \delta^{-1}((\delta y)^k - 1)$. If $a(0)\mathrm{sgn}\,(\tau) \notin \{1, 3, 5 \dots\}$ then the analysis of section 3 establishes that the absolute value of any eigenvalue of $|w|^{-2\sigma} D_{w,0}$ is bounded below, uniformly as $\delta \to 0^+$ (equivalently, as $w \to -\infty$).

If $a(0)\mathrm{sgn}\,(\tau) \in \{3, 5, \dots\}$ then the situation does degenerate as $\delta \to 0$, but for all $k > 1$ Lemma 4.1 guarantees that all eigenvalues of $|w|^{-2\sigma} D_{w,0}$ are uniformly bounded below by $c\delta^2$ as $\delta \to 0$, for some $c > 0$. □

For any $m, n$, $\partial^m \lambda_n(w, 0)/\partial \epsilon^m$ is also an entire holomorphic function of $w$, since it is locally a uniform limit of iterated difference quotients of the holomorphic functions $w \mapsto \lambda_n(w, \epsilon)$.

**Lemma 5.3.** *If $k > 1$, $\operatorname{sgn}(\tau)a(0) = +1$ and $a^{(m)}(0) \ne 0$ for some $m \ge 1$ then for each $n$ there exists $0 \le \nu \le m$ such that the function $w \mapsto \partial^\nu \lambda_n/\partial \epsilon^\nu(w,0)$ does not vanish identically.*

*Proof.* Setting $y = \delta^{-1} + t$,
$$|w|^{-2\sigma} D_{w,\epsilon} = -\partial_t^2 + q_\delta^2(t) - q_\delta'(t) - b(\epsilon \delta^{r-1} + \epsilon \delta^r t)q_\delta'(t)$$
where $r = \sigma/(1-\sigma) \in (0,1)$. If we restrict attention to the case where not only $\epsilon, \delta$ but also the larger quantity $\epsilon \delta^{r-1}$ remain in a sufficiently small neighborhood of the origin then the analysis of section 3 establishes that $D_{w,\epsilon}$ has a unique small eigenvalue $\tilde\lambda(w,\epsilon)$, which takes the form $\tilde\lambda(w,\epsilon) = h(z,\delta)$ where $h \in C^\infty$ near 0 is the unique small eigenvalue of $B_{z,\delta}$ with $z = \epsilon \delta^{r-1}$. By Corollary 3.10, if $a^{(m)}(0) \ne 0$, then $\partial^m h/\partial z^m \ne 0$ in a neighborhood of the origin. Since $\tilde\lambda(w,\epsilon) = h(\epsilon \delta^{r-1}, \delta)$ and $\delta$ is a function of $w$ alone,
$$\partial^m \tilde\lambda/\partial \epsilon^m = \delta^{m(r-1)} \partial^m h/\partial z^m.$$
This is nonzero wherever $\epsilon \delta^{r-1}, \delta$ are sufficiently small; in particular, is nonzero at $(w,\epsilon)$ whenever $\epsilon = 0$ and $|w| = \delta^{-1/(1-\sigma)}$ is sufficiently large.

So far we have treated only one eigenvalue. But since the eigenvalues are distinct and all vary holomorphically, the conclusion of the lemma holds automatically with $\nu = 0$ for all except at most one index $n$, namely that index corresponding to the unique small eigenvalue (of $D_{w,0}$) for the range of $w$ just discussed. □

Consider lastly the special and simplest case $k = 1$. Then substituting $x = y - w$,
$$D_{w,\epsilon} = -\partial_x^2 + (x+w)^2 - \operatorname{sgn}(\tau)[a(0) + b(\epsilon x)]$$
$$= -\partial_y^2 + y^2 - \operatorname{sgn}(\tau)[a(0) + b(\epsilon y - \epsilon w)].$$
Again each eigenvalue is an entire holomorphic function of $w$, uniformly for $\epsilon$ in any compact set. If $a(0)\operatorname{sgn}(\tau) \notin \{1,3,5,\dots\}$ then as for the case $k > 1$, there clearly exists $\delta > 0$ such that $|\lambda_n(w,\epsilon)| \ge \delta$ for all sufficiently small $(w,\epsilon)$ and all $n$.

If $\operatorname{sgn}(\tau)a(0) \in \{1,3,5,\dots\}$ write $D_{w,\epsilon} = -\partial_y^2 + y^2 - \operatorname{sgn}(\tau)[a(0) + b(\epsilon y - z)]$, with $z = \epsilon w$. Let $\lambda(z,\epsilon)$ be the unique small eigenvalue of this last operator, for all sufficiently small $(z,\epsilon)$. Then $\lambda(\epsilon w, \epsilon)$ is the unique small eigenvalue $\tilde\lambda(w,\epsilon)$ of $D_{w,\epsilon}$ for small $\epsilon$ and bounded $w$.

**Lemma 5.4.** *If $k = 1$, $a(0)\operatorname{sgn}(\tau) \in \{1,3,5,\dots\}$ and $a^{(m)}(0) \ne 0$ for some $m \ge 1$, then*
$$\frac{\partial^m \lambda}{\partial^m z}(0,0) \ne 0.$$

*Proof.* $\lambda(z,0) \equiv -\operatorname{sgn}(\tau) b(-z) = \pm[a(-z) - a(0)]$. □

**Corollary 5.5.** *If $k = 1$, $a(0)\mathrm{sgn}\,(\tau) \in \{1, 3, 5, \ldots\}$ and $a^{(m)}(0) \neq 0$ for some $m \geq 1$, then there exists $n \geq 0$ such that the function $w \mapsto \partial^n \tilde{\lambda}/\partial \epsilon^n(w, 0)$ is nonconstant.*

*Proof.* Fix any $n \geq 0$ for which $\partial^n \lambda / \partial z^n(0, 0) \neq 0$.

$$\frac{\partial^n \tilde{\lambda}}{\partial w^n}(w, \epsilon) = \epsilon^n \frac{\partial^n \lambda}{\partial z^n}(\epsilon w, \epsilon) = c\epsilon^n + O(\epsilon^{n+1})$$

for some $c \neq 0$. Thus

$$\frac{\partial^n}{\partial w^n}\frac{\partial^n \tilde{\lambda}}{\partial \epsilon^n}(0, 0) \neq 0. \quad \square$$

**Lemma 5.6.** *Suppose that $a^{(m)}(0) \neq 0$ for some $m \geq 1$, or that $a(0)$ does not belong to $\{\pm 1, \pm 3, \pm 5 \ldots\}$, or that $k > 1$ and that $a(0) \notin \{\pm 1\}$. Then there exists $M < \infty$ such that*

$$\left|\{(\eta, \tau) \in \mathcal{D}_k : \min_j |\mu_j(\eta, \tau)| \leq |\tau|^{-M}\}\right| < \infty.$$

As always, $\{\mu_j\}$ denote the eigenvalues of $A_{\eta,\tau}$.

*Proof.* Consider $S_q = \{(\eta, \tau) \in \mathcal{D}_k : 2^q \leq |\tau| \leq 2^{q+1}\}$ for each nonnegative integer $q$. The Jacobian determinant for the change of variables $(\eta, \tau) \mapsto (w, \epsilon)$ is

$$\left|\frac{\partial(\eta, \tau)}{\partial(w, \epsilon)}\right| = C|\tau|^{(k+2)(k+1)}|\tau|^{1/(k+1)} \leq C 2^{2q}.$$

Set $\Omega = \{(w, \epsilon) : |w| \leq \gamma_0 \text{ and } 0 \leq \epsilon \leq \tau_0^{-1/(k+1)}\}$. It suffices to show that

$$\sum_q 2^{2q}|\{(w, \epsilon) \in \Omega : \min_j |\lambda_j(w, \epsilon) \leq 2^{-Mq}\}| < \infty,$$

provided that $\tau_0$ and $M$ are chosen to be sufficiently large.

Fix any $j$. Since any nonconstant analytic function has some nonvanishing derivative at each point, the preceding lemmas guarantee that for each point $(w, 0) \in \Omega$ there exists some multi-index $\alpha$ (possibly equal to $(0, 0)$) such that

$$[\partial^\alpha \lambda_j / \partial(w, \epsilon)^\alpha](w, 0)$$

is nonzero. If $\tau_0$ is chosen to be sufficiently large, then the same holds at each $(w, \epsilon) \in \Omega$, since $w$ ranges only over a compact set. $\Omega$ may then be partitioned into finitely many two-dimensional rectangles $\Omega_i$, in each of which some partial derivative $\partial^\alpha \lambda_j / \partial(w, \epsilon)^\alpha$ is nonzero, with $\alpha$ depending on $i$ but not otherwise on $(w, \epsilon)$. Lemma 3.11 then implies a lower bound

$$|\{(w, \epsilon) \in \Omega_i : |\lambda_j(w, \epsilon)| \leq 2^{-Mq}\}| \leq C 2^{-\delta M q}$$

for some $\delta > 0$. Choosing $M$ to be sufficiently large relative to $\delta$,
$$\sum_q 2^{2q} C 2^{-\delta M q} < \infty,$$
and the proof would be complete if we were interested only in one eigenvalue $\lambda_j$ rather than in their minimum.

Recall that each $\lambda_j$ is a continuous function on $\Omega$, and that for any fixed $\omega \in \Omega$, $\lambda_j(\omega) \to +\infty$ as $j \to \infty$. Fix $N(\omega)$ such that $\lambda_j(\omega) > 1$ for all $j \geq N(\omega)$. Since $\lambda_0 < \lambda_1 < \ldots$ at every point, there exists some neighborhood $V$ of $\omega$ such that $\lambda_j(w, \epsilon) \geq 1$ for every $(w, \epsilon) \in V$, for every $j \geq N(\omega)$. Since $\Omega$ is compact, there exists $N' < \infty$ such that $\lambda_j(w, \epsilon) \geq 1$ for every $j > N'$, for every $(w, \epsilon) \in \Omega$. Thus only finitely many eigenvalues $\lambda_0, \ldots \lambda_{N'}$ need be taken into account in analyzing the minimum (in absolute value) eigenvalue, so the result follows from the preceding paragraph. $\square$

## 6. Nonsolvable cases

**Proposition 6.1.** *If $k > 1$, $a(0) \in \{\pm 1\}$ and $a^{(m)}(0) = 0$ for all $m \geq 1$ then $L$ is not locally solvable at $0$.*

*Proof.* Throughout the discussion we assume that $a(0) = +1$; the case $a(0) = -1$ reduces to this by the change of variables $(x, y, t) \mapsto (x, -y, -t)$. In all cases we replace $y$ by $-y$, thus converting $\partial_y$ to $-\partial_y$. For $x$ near $0$ we are then dealing with a small perturbation of
$$L_0 = -\partial_x^2 - (-\partial_y + x^k \partial_t)^2 + ia(0)kx^{k-1}\partial_t$$
$$= (-\partial_x - i(-\partial_y + x^k \partial_t))(\partial_x - i(-\partial_y + x^k \partial_t)).$$

Throughout this proof it is assumed that $\eta, \tau$ are both positive. Define
$$g_{\eta,\tau}(x) = \exp(\eta x - \tau(k+1)^{-1} x^{k+1}).$$
Then
$$L_0 \left( e^{i\eta y + i\tau t} g_{\eta,\tau} \right) \equiv 0. \tag{6.1}$$

$g_{\eta,\tau}$ is a Schwartz function for odd $k$, but not so for even $k$, and this will complicate the formulas to follow. $g_{\eta,\tau}$ has a critical point at $x = (\eta/\tau)^{1/k}$, where we take the unique positive root. The critical value is
$$g_{\eta,\tau}((\eta/\tau)^{1/k}) = \exp\left( \eta^{(k+1)/k} \tau^{-1/k} - (k+1)^{-1} \eta^{(k+1)/k} \tau^{-1/k} \right)$$
$$= \exp\left( \frac{k}{k+1} \eta^{(k+1)/k} \tau^{-1/k} \right),$$

so we normalize by setting
$$G_{\eta,\tau}(x) = \exp\left(\frac{-k}{k+1}\eta^{(k+1)/k}\tau^{-1/k}\right)g_{\eta,\tau}(x)$$

so that $G_{\eta,\tau}((\eta/\tau)^{1/k}) \equiv 1$.

Fix a cutoff function $h \in C_0^\infty(\mathbf{R})$ satisfying $h(0) = 1$, everywhere nonnegative and supported in $[-1/2, 1/2]$. Let $\lambda \in \mathbf{R}^+$ be a large parameter, which will eventually be allowed to tend to $\infty$. Define
$$F_\lambda(x,y,t) = \iint_{\mathbf{R}^2} e^{i(\eta y + \tau t)}G_{\eta,\tau}(x)h(\lambda^{-3/4}(\tau - \lambda))h(\lambda^{-1/4}(\eta - \lambda^{1/2}))\,d\eta\,d\tau.$$

From (6.1) it follows that for all $\lambda$,
$$L_0 F_\lambda \equiv 0.$$

In the support of the integrand $\eta, \tau$ satisfy
$$|\tau - \lambda| \leq \lambda^{3/4} \text{ and } |\eta - \lambda^{1/2}| \leq \lambda^{1/4}. \tag{6.2}$$

Throughout the remainder of this proof $\eta, \tau$ are always assumed to satisfy these last two inequalities[4].

The critical point $(\eta/\tau)^{1/k}$ satisfies
$$|(\eta/\tau)^{1/k} - \lambda^{-1/2k}| \leq C\lambda^{-1/4}\lambda^{-1/2k}.$$

If
$$\frac{1}{6}(\eta/\tau)^{1/k} \leq |x - (\eta/\tau)^{1/k}| \leq \frac{5}{6}(\eta/\tau)^{1/k}$$

then
$$|G_{\eta,\tau}| \leq \exp(-c\lambda^{\frac{1}{2}\frac{k+1}{k} - \frac{1}{k}}) \leq \exp(-c\lambda^c)$$

for some $c > 0$, uniformly in $\eta, \tau, \lambda$. For all sufficiently large $\lambda$, this holds for all $x \in I_\lambda$ where
$$I_\lambda = \{x : \frac{1}{3}\lambda^{-1/2k} \leq |x - x_\lambda| \leq \frac{2}{3}\lambda^{-1/2k}\} \tag{6.3}$$

and
$$x_\lambda = \lambda^{-1/2k}.$$

More generally, for any multi-index $\alpha$,
$$\left|\frac{\partial^\alpha G_{\eta,\tau}(x)}{\partial(x,\eta,\tau)^\alpha}\right| \leq \exp(-c\lambda^c)$$

---

[4]These precise exponents $1/4, 2/4, 3/4, 1$ have been chosen for algebraic simplicity and have no intrinsic significance

for all $x \in I_\lambda$, uniformly in $\eta, \tau$ satisfying (6.2), for some $c > 0$ depending on $\alpha$, for all sufficiently large $\lambda$. The same conclusion then follows with $G_{\eta,\tau}$ replaced by $\exp(i\eta y + i\tau t) G_{\eta,\tau}(x)$, by Leibniz's rule. For all $|x - x_\lambda| \leq \frac{2}{3}\lambda^{-1/2k}$ there is the weaker bound

$$\left|\frac{\partial^\alpha G_{\eta,\tau}(x)}{\partial(x, \eta, \tau)^\alpha}\right| \leq \lambda^{C_\alpha}$$

for some $C_\alpha < \infty$.

These upper bounds for $G_{\eta,\tau}$ and its partial derivatives lead to corresponding bounds for $F_\lambda$. For $|x - x_\lambda| \leq \frac{2}{3}\lambda^{-1/2k}$,

$$\left|\frac{\partial^{\alpha+\beta+\gamma} F_\lambda}{\partial x^\alpha \partial y^\beta \partial t^\gamma}\right| \leq \iint_{\substack{|\tau-\lambda| \leq \lambda^{3/4} \\ |\eta-\lambda^{1/2}| \leq \lambda^{1/4}}} |\partial_x^\alpha G_{\eta,\tau}(x)| \eta^\beta \tau^\gamma \, d\eta \, d\tau \leq \lambda^C$$

for large $\lambda$, where $C$ depends only on $\alpha, \beta, \gamma$. For $x \in I_\lambda$ there is the improved bound

$$\left|\frac{\partial^{\alpha+\beta+\gamma} F_\lambda}{\partial x^\alpha \partial y^\beta \partial t^\gamma}\right| \leq e^{-\lambda^c}$$

for some $c > 0$.

$F_\lambda$ is also very small, for large $\lambda$, if $(y, t)$ is not very close to the origin. More precisely, integrating by parts $N$ times with respect to $\tau$ in the integral defining $F_\lambda$ gives, for $|x - x_\lambda| \leq \frac{2}{3}\lambda^{-1/2k}$,

$$|F_\lambda(x, y, t)| \leq C_N |t|^{-N} \iint_{\substack{|\tau-\lambda| \leq \lambda^{3/4} \\ |\eta-\lambda^{1/2}| \leq \lambda^{1/4}}} G_{\eta,\tau}(x) \, d\eta \, d\tau \leq C_N \lambda^{C_0} \lambda^{-N/2} |t|^{-N}$$

with $C_0$ independent of $N$. Indeed, consider

$$\frac{\partial}{\partial \tau}\left(G_{\eta,\tau}(x) h(\lambda^{-3/4}(\tau - \lambda)) h(\lambda^{-1/4}(\eta - \lambda^{1/2}))\right).$$

When the derivative falls on the normalizing factor $\exp(-\eta^{(k+1)/k}\tau^{-1/k}k/(k+1))$, the result is an additional factor of $(\eta^{(k+1)/k}\tau^{-1-1/k}) = O(\lambda^{-(k+1)/2k})$. When it falls on $g_{\eta,\tau}(x)$, the result is a factor of $(k+1)^{-1}x^{k+1} = O(\lambda^{-(k+1)/2k})$. When it falls on $h(\lambda^{-3/4}(\tau - \lambda))$, the result is $O(\lambda^{-3/4})$. A second derivative with respect to $\tau$ either falls again on $G_{\eta,\tau}(x) h(\lambda^{-3/4}(\tau - \lambda)) h(\lambda^{-1/4}(\eta - \lambda^{1/2}))$, producing a second factor that is $O(\lambda^{-(k+1)/2k})$, or falls on the factor $\eta^{(k+1)/k}\tau^{-1-1/k}$, netting another factor of $\tau^{-1} \sim \lambda^{-1}$. Thus each derivative nets a factor smaller than a constant times $\lambda^{-1/2}$. Iterating $N$ times, we obtain a bound of $C_N(\lambda^{1/2}|t|)^{-N}$.

Integrating by parts instead $N$ times with respect to $\eta$ and applying the same reasoning gives

$$|F_\lambda(x, y, t)| \leq C_N \lambda^{C_0} \lambda^{-N/2k} |y|^{-N}.$$

The same bounds hold for $\partial^\alpha F_\lambda / \partial(x,y,t)^\alpha$ with an extra factor of $C_{N,\alpha} \lambda^{C_\alpha}$ for each $\alpha$, while for $x \in I_\lambda$ there is an additional factor of $\exp(-\lambda^c)$ for some $c > 0$. Our primary conclusion is then that for each $\alpha$,

$$\frac{\partial^\alpha F_\lambda}{\partial(x,y,t)^\alpha} = O(\lambda^{-M}) \tag{6.4}$$

for all $M < \infty$ uniformly in $(x,y,t), \lambda$ where $x \in I_\lambda$ or $|t| \geq \lambda^{-1/4}$ or $|y| \geq \lambda^{-1/4k}$.

A crude lower bound on $F_\lambda$ will also be required. If $|\tau - \lambda| \leq \lambda^{-1}$ and $|\eta - \lambda^{1/2}| \leq \lambda^{-1}$ then $\tau^{-1/k} = \lambda^{-1/k} + O(\lambda^{-2})$ and $\eta^{1/k} = \lambda^{1/2k} + O(\lambda^{-3/2})$, so $x_\lambda - (\eta/\tau)^{1/k} = O(\lambda^{-3/2})$. Consequently

$$\left| \log \frac{g_{\eta,\tau}(x_\lambda)}{g_{\eta,\tau}((\eta/\tau)^{1/k})} \right| \leq C\lambda^{1/2} \lambda^{-3/2} + C\lambda(\lambda^{-k/2k} \lambda^{-3/2}) = O(1),$$

so $G_{\eta,\tau}(x_\lambda)$ is bounded below by a strictly positive constant independent of $\lambda$, for all such $(\eta, \tau)$. Thus there exists $c > 0$ such that

$$F_\lambda(x_\lambda, 0, 0) \geq c \iint_{\substack{|\tau - \lambda| \leq \lambda^{-1} \\ |\eta - \lambda^{1/2}| \leq \lambda^{-1}}} 1 \, d\eta \, d\tau \geq c\lambda^{-2}.$$

Since $\nabla F_\lambda = O(\lambda^C)$ for some finite $C$, there consequently exists $B \in \mathbf{R}^+$ such that for all sufficiently large $\lambda$,

$$F_\lambda(x,y,t) \geq c'\lambda^{-2} \text{ whenever } |(x,y,t) - (x_\lambda, 0, 0)| \leq \lambda^{-B}. \tag{6.5}$$

A necessary condition [7] for any linear operator $\mathcal{L}$ to be solvable at 0 is that there exist $\epsilon > 0$, $N < \infty$ such that for all $\phi, \psi \in C_0^\infty(\mathbf{R}^3)$ supported in $\{|(x,y,t)| \leq \epsilon\}$,

$$\left| \int \phi \psi \right| \leq N \|\phi\|_{C^N} \|\mathcal{L}^* \psi\|_{C^N} \tag{6.6}$$

where $\mathcal{L}^*$ denotes the transpose of $\mathcal{L}$. We will prove that (6.6) does not hold for $\mathcal{L} = L^*$; thus $L^*$ is not locally solvable. Since the class of operators under discussion in Proposition 6.1 is closed under taking transposes, this will conclude the proof.

Fix a cutoff function $\zeta \in C_0^\infty(\mathbf{R})$ supported in $[-2/3, 2/3]$, such that $\zeta(s) \equiv 1$ for $|s| \leq 1/3$. For large $\lambda$ set

$$\psi_\lambda(x,y,t) = F_\lambda(x,y,t) \zeta_\lambda(x,y,t)$$

where

$$\zeta_\lambda(x,y,t) = \zeta((x - x_\lambda)/\lambda^{-1/2k}) \zeta(|t|/\lambda^{-1/8}) \zeta(|y|/\lambda^{-1/8k}).$$

The gradient of $\zeta_\lambda$ is supported in a region where $\partial^\alpha F_\lambda / \partial(x,y,t)^\alpha = O(\lambda^{-M})$ for every finite $M$ and every $\alpha$, by (6.4).

Fix any finite exponent $N$. Choose a nonnegative test function $\phi_\lambda \in C_0^\infty$ supported where $|(x,y,t) - (x_\lambda, 0, 0)| \leq \lambda^{-B}$, with $\phi_\lambda(x,y,t) \equiv 1$ where $|(x,y,t) - (x_\lambda, 0, 0)| \leq \frac{1}{2}\lambda^{-B}$, satisfying $\|\phi_\lambda\|_{C^N} = O(\lambda^{NB})$. (6.5) thus implies

$$\int \phi_\lambda \psi_\lambda \geq \delta \lambda^{-3B-2}$$

for some $\delta > 0$.

In order to prove that $L^*$ is not locally solvable at $0$ we aim to prove that (6.6) is violated, for the arbitrary exponent $N$ already introduced, for all sufficiently large $\lambda$. To do this it now suffices to prove

$$\|L\psi_\lambda\|_{C^N} = O(\lambda^{-A}) \text{ for all } A < \infty.$$

Recall that $b(x) = a(x) - a(0)$ and $L = L_0 + ikx^{k-1}b(x)\partial_t$. Then denoting by $\zeta_\lambda$ also the operator defined by multiplication by the function $\zeta_\lambda$, and recalling that $L_0 F_\lambda \equiv 0$,

$$L\psi_\lambda = L_0(F_\lambda \zeta_\lambda) + ikx^{k-1}b(x)\partial_t(F_\lambda \zeta_\lambda)$$
$$= [L_0, \zeta_\lambda]F_\lambda + O(|b(x)| \cdot |\nabla(\zeta_\lambda F_\lambda)|)$$

uniformly at all points of the support of $\psi_\lambda$. $\zeta_\lambda F_\lambda$ is supported where $|x| \leq C\lambda^{-1/2k}$ and is $O(\lambda^C)$ in $C^1$ norm for some finite $C$. Since $|b(x)| = O(|x|^R)$ for all $R < \infty$, the final term in the last display is $O(\lambda^{-A})$ for all $A < \infty$. The differential operator $[L_0, \zeta_\lambda]$ is of order one, and has smooth coefficients supported in the union of the three regions where $x \in I_\lambda$ or $|t| \geq \lambda^{-1/8}$ or $|y| \geq \lambda^{-1/8k}$. In supremum norm these coefficients are $O(\lambda)$. (6.4) therefore guarantees that $[L_0, \zeta_\lambda]F_\lambda$ is likewise $O(\lambda^{-A})$ in the $C^0$ norm, for all finite exponents $A$.

The same reasoning applies to the $C^N$ norm, for any finite $N$. This completes the proof of Proposition 6.1. □

**Proposition 6.2.** *If $k = 1$, $a(0) \in \{\pm 1, \pm 3, \pm 5 \ldots\}$ and $a^{(m)}(0) = 0$ for all $m \geq 1$ then $L$ is not locally solvable at $0$.*

*Proof.* Write $L = L_0 + ib(x)\partial_t$. A much simpler version of the above reasoning shows that there exists a Schwartz function $F$ in $\mathbf{R}^3$ satisfying $L_0 F \equiv 0$ and $F(0) \neq 0$. Setting $F_\lambda(x, y, t) = F(\lambda x, \lambda y, \lambda^2 t)$, $L_0 F_\lambda \equiv 0$ for all $\lambda \in \mathbf{R}^+$. Define now $\psi_\lambda(x, y, t) = F_\lambda(x, y, t)\zeta(\lambda^{1/2}x)\zeta(\lambda^{1/2}y)\zeta(\lambda t)$. Since $F_\lambda$ belongs to the Schwartz class and $L_0 F_\lambda \equiv 0$, $L\psi_\lambda = O(\lambda^{-A})$ in the $C^N$ norm, for any $N, A < \infty$. Define $\phi_\lambda$ to be $\phi(\lambda x, \lambda y, \lambda^2 t)$ for some fixed nonnegative $\phi \in C_0^\infty(\mathbf{R}^3)$ that is supported in a sufficiently small neighborhood of the origin and satisfies $\int \phi \neq 0$. Then (6.6) is violated by this pair $\psi_\lambda, \phi_\lambda$ for all sufficiently large $\lambda$, for any given $N$. □

As mentioned in Corollary 1.2, our theory does include locally solvable operators that are not locally solvable in $L^2$.

**Remark.** Suppose that $k = 1$ and $a(0) \in \{\pm 1, \pm 3, \pm 5, \dots\}$, or that $k > 1$ is odd and $a(0) \in \{\pm 1\}$. If $a^{(1)}(0) = a^{(2)}(0) = 0$, then $L$ is not locally solvable in $L^2$ at 0.

*Proof.* In these cases the basic operator $L_0 = -X^2 - Y^2 + ia(0)[X, Y]$ has the property that there exists a function $f$ not vanishing identically, belonging to the Schwartz class in $\mathbf{R}^3$, and satisfying $L_0 f \equiv 0$. Indeed, either for all $\tau > 0$, or for all $\tau < 0$, the ordinary differential operator $A_{\eta,\tau}$ corresponding to $L_0$ annihilates a function $f_{\eta,\tau}$ in the Schwartz class on $\mathbf{R}^1$, for all $\eta$. $f$ is then defined as the inverse partial Fourier transform of $h(\eta, \tau) f_{\eta,\tau}(x)$, for some $h \in C_0^\infty(\mathbf{R}^2)$.

Fix a cutoff function $\zeta \in C_0^\infty(\mathbf{R}^3)$ that is identically equal to 1 in some neighborhood of the origin, and define
$$F_\lambda(x, y, t) = \lambda^{(k+3)/2} \cdot f(\lambda x, \lambda y, \lambda^{k+1} t) \zeta(x, y, t)$$
for each large $\lambda \in \mathbf{R}^+$. Then $\|F_\lambda\|_{L^2}$ equals a constant modulo $O(\lambda^{-N})$ for all $N$. Clearly $\|L_0 F_\lambda\|_{L^2} = O(\lambda^{-N})$ for all $N$, since $F_\lambda$ and all of its derivatives are $O(\lambda^{-N})$ on the support of $\nabla \zeta$. The $L^2$ norm of $\partial_t F_\lambda$ is $O(\lambda^{k+1})$, and $F_\lambda$ is essentially supported where $x = O(\lambda^{-1})$, so
$$\|[a(x) - a(0)]x^{k-1} \partial_t F_\lambda\|_{L^2} \leq C\lambda^{-3}\lambda^{-(k-1)}\lambda^{k+1},$$
assuming that $a(x) - a(0) = O(x^3)$. In all, $\|LF_\lambda\| = O(\lambda^{-1})$ as $\lambda \to \infty$, so $\|F_\lambda\| \gg \|L^* F_\lambda\|$, since $L$ equals its transpose. Because $\zeta$ may be taken to be supported in any given neighborhood of 0, by duality this implies local nonsolvability in $L^2$. □

University of California, Los Angeles, Los Angeles, Ca. 90095-1555
*E-mail address*: `christ@math.ucla.edu`

Institute of Mathematics, Bulgarian Academy of Science, P.O. Box 373, 1090 Sofia, Bulgaria
*E-mail address*: `difeq@bgearn.bitnet`